\RequirePackage{rotating}

\documentclass{amsart}%
\usepackage{amsfonts}
\usepackage{amsmath}
\usepackage{amssymb}
\usepackage{amsthm}
\usepackage{booktabs}
\usepackage{mathtools}
\usepackage[colorlinks, linkcolor=red]{hyperref}
\usepackage[noabbrev,nameinlink]{cleveref}
\usepackage{enumerate}
\usepackage{bm}
\usepackage{tikz}
\usepackage{tikz-cd}
\usetikzlibrary{decorations.pathreplacing,patterns}
\usepackage{subcaption}
\usepackage{xcolor}
\usepackage{rotating}
\usepackage[boxsize=1em]{ytableau}

\makeatletter
\newtheorem*{rep@theorem}{\rep@title}
\newcommand{\newreptheorem}[2]{%
\newenvironment{rep#1}[1]{%
 \def\rep@title{#2 \ref{##1}}%
 \begin{rep@theorem}}%
 {\end{rep@theorem}}}
\makeatother

\makeatletter
\newtheorem*{rep@corollary}{\rep@title}
\newcommand{\newrepcorollary}[2]{%
\newenvironment{rep#1}[1]{%
 \def\rep@title{#2 \ref{##1}}%
 \begin{rep@corollary}}%
 {\end{rep@corollary}}}
\makeatother

\newtheorem{theorem}{Theorem}[section]
\newreptheorem{theorem}{Theorem}

\newtheorem{corollary}[theorem]{Corollary}
\newtheorem{lemma}[theorem]{Lemma}
\newreptheorem{corollary}{Corollary}
\newtheorem{problem}[theorem]{Problem}

\theoremstyle{definition}
\newtheorem{definition}[theorem]{Definition}
\newtheorem{notation}[theorem]{Notation}
\newtheorem{remark}[theorem]{Remark}
\newtheorem{example}[theorem]{Example}
\newtheorem{question}[theorem]{Question}

\theoremstyle{plain}
\newtheorem*{namedthm}{\namedthmname}
\newcounter{namedthm}



\newcommand{\bC}{\mathbb{C}}

\newcommand{\bQ}{\mathbb{Q}}

\newcommand{\bZ}{\mathbb{Z}}

\newcommand{\cH}{\mathcal{H}}
\newcommand{\cI}{\mathcal{I}}
\newcommand{\cJ}{\mathcal{J}}

\newcommand{\cR}{\mathcal{R}}

\newcommand{\cSDR}{\mathcal{SDR}}
\newcommand{\cSH}{\mathcal{SH}}
\newcommand{\cSR}{\mathcal{SR}}
\newcommand{\cX}{\mathcal{X}}

\newcommand{\fS}{\mathfrak{S}}

\DeclareMathOperator{\Ann}{Ann}

\DeclareMathOperator{\coinv}{coinv}
\newcommand{\dif}{\mathrm{d}}

\DeclareMathOperator{\GrFrob}{GrFrob}

\DeclareMathOperator{\id}{id}

\DeclareMathOperator{\inv}{inv}

\DeclareMathOperator{\maj}{maj}

\DeclareMathOperator{\rev}{rev}

\DeclareMathOperator{\sgn}{sgn}
\DeclareMathOperator{\Span}{Span}
\DeclareMathOperator{\Stir}{Stir}

\DeclareMathOperator{\Sum}{sum}

\DeclareMathOperator{\wgt}{wgt}

\newlength{\cellsize} 
\newsavebox{\cell}
\sbox{\cell}{\begin{picture}(18,18)
\put(0,0){\line(1,0){18}}
\put(0,0){\line(0,1){18}}
\put(18,0){\line(0,1){18}}
\put(0,18){\line(1,0){18}}
\end{picture}}
\newsavebox{\dashedcell}
\sbox{\dashedcell}{\begin{picture}(18,18)
\multiput(0,0)(5,0){4}{\line(1,0){3}}
\multiput(0,0)(0,5){4}{\line(0,1){3}}
\multiput(18,0)(0,5){4}{\line(0,1){3}}
\multiput(0,18)(5,0){4}{\line(1,0){3}}
\end{picture}}
\newcommand\cellify[1]{\def\thearg{#1}\def\nothing{}%
\ifx\thearg\nothing
\hbox to 0pt{ \hss}\else%
\hbox to 0pt{\usebox{\cell} \hss}\fi%
\vbox to \cellsize{
\vss
\hbox to \cellsize{\hss$#1$\hss}
\vss}}
\newcommand\tableau[1]{\vtop{\let\\\cr
\baselineskip -16000pt \lineskiplimit 16000pt \lineskip 0pt
\ialign{&\cellify{##}\cr#1\crcr}}}

\title[Tanisaki witness relations]{Tanisaki witness relations \\ for harmonic differential forms}
\author{Joshua P. Swanson}
\date{\today}

\begin{document}
  
\begin{abstract}
  Inspired by a series of conjectures and formulas related to higher coinvariant algebras, we present two families of relations involving harmonic differential forms of the symmetric group. Our relations, together with a novel bijection, are sufficient to give a filtration of the $1$-forms suggested by work of Haglund--Rhoades--Shimozono with composition factors given by Tanisaki quotients. These are ``almost all'' of the necessary relations in a certain asymptotic sense we make precise.

\end{abstract}
\maketitle

\setcounter{tocdepth}{1}
\tableofcontents

\section{Introduction}

\subsection{Overview of results}

We present a large family of relations between harmonic differential forms of the symmetric group. These relations are involved in a series of conjectures and results concerning higher coinvariant algebras. Combining a recent conjecture of Zabrocki \cite{1902.08966} on super diagonal coinvariant algebras with results of Haglund--Rhoades--Shimozono \cite{MR3783430} on generalized coinvariant algebras related to the Delta Conjecture of Haglund--Remmel--Wilson \cite{MR3811519} suggests the existence of a filtration of the harmonic differential forms of the symmetric group whose successive quotients are cohomology rings of Springer fibers; see \Cref{que:main}. These cohomology rings were given a well-known presentation by Tanisaki \cite{MR685425}. Separately, a recent conjecture of Wallach and the author \cite{sw2}\footnote{\label{foot:FW}This conjecture was proven by Rhoades--Wilson \cite{RW23} after this manuscript was submitted.} gives an explicit description of the harmonic differential forms in terms of certain differential operators from \cite{MR4258767} applied to the Vandermonde determinant.

Our main results are two families of ``Tanisaki witness relations'' between these explicit harmonic differential forms, \Cref{thm:generic_Pieri} and \Cref{thm:hook_witness_extreme}. Together with a novel bijection, \Cref{thm:bijection}, our relations are sufficient to prove the hoped-for filtration for $1$-forms, \Cref{cor:1_forms}. They also provide ``almost all'' of the necessary relations in a certain asymptotic sense, see \Cref{rem:generic}. Our arguments are combinatorial and effectively construct certain intricate sign-reversing involutions.

Our results provide further evidence for the above conjectures and precisely identify some of the remarkably rich structure underlying them. We hope they will spur additional research on this topic, especially from topological, homological, algebraic, or geometric perspectives.

The rest of this introduction describes these developments and their context in detail and states our main results. We introduce the classical coinvariants, the Tanisaki ideals, the generalized coinvariant algebras, the Delta Conjecture and higher coinvariant algebras, classical harmonics, harmonic differential forms, the flip action, the potential filtration which motivated this work, Tanisaki witness relations, our first family of relations which we call the Generic Pieri Rule, our result for $1$-forms, and our second family of relations.

\subsection*{Acknowledgements}

\Cref{que:main}, the observation \eqref{eq:Phi_condition}, and the first computationally discovered Tanisaki witness relations are due to Brendon Rhoades, without whom this paper would not exist and whom I heartily thank. Thanks also go to Nolan Wallach for collaboration on related projects and Dani\"el Kroes for useful discussion on the bijection $\Phi_n$ from \Cref{sec:bijection}. Finally, I must also thank the anonymous referee for their careful reading of the manuscript and excellent suggestions.

\subsection{Classical coinvariants}

The \textit{coinvariant algebra} of the symmetric group $\fS_n$ is the quotient
\begin{equation}
  \cR_n \coloneqq \frac{\mathbb{Q}[x_1, \ldots, x_n]}{\langle e_r(\underline{n}) : r \in [n]\rangle},
\end{equation}
where the generators are the \textit{elementary symmetric polynomials}
  \[ e_r(\underline{n}) \coloneqq e_r(x_1, \ldots, x_n) \coloneqq \sum_{1 \leq i_1 < \cdots < i_r \leq n} x_{i_1} \cdots x_{i_r}. \]
The classical coinvariant algebra is very well understood from topological, geometric, and combinatorial perspectives \cite{MR2538310,MR1464693,MR597728,MR4029695,MR526968}. As one example, Borel \cite{MR51508} showed that $\cR_n$ is a presentation for the cohomology ring of the \textit{complete flag variety}. 

\subsection{Tanisaki ideals}

A series of authors (see \cite[p.83]{MR1168926}) more generally considered the cohomology ring of the \textit{Springer fiber} $\cX_\mu$ consisting of complete flags in $\mathbb{C}^n$ fixed by the unipotent matrix with Jordan blocks of size $\mu_1 \geq \cdots \geq \mu_n \geq 0$ for a partition $\mu \vdash n$. The complete flag variety is the case $\mu=(1, \ldots, 1)$, corresponding to the identity matrix.

Tanisaki \cite{MR685425} gave a presentation of the cohomology rings $H^*(\overline{\cX_\mu})$,
\begin{equation}
  \cR_\mu \coloneqq \frac{\bQ[x_1, \ldots, x_n]}{\cI_\mu}
\end{equation}
where
\begin{equation}
  \cI_\mu \coloneqq \langle e_r(S) : |S| - d_{|S|}(\mu) < r \leq |S|, S \subset [n]\rangle
\end{equation}
is a \textit{Tanisaki ideal}, with $d_k(\mu) \coloneqq \mu_n' + \mu_{n-1}' + \cdots + \mu_{n-k+1}'$ and
\begin{equation}\label{eq:Tanisaki_generators}
  e_r(S) \coloneqq \sum_{\{i_1 < \cdots < i_r\} \subset S} x_{i_1} \cdots x_{i_r}.
\end{equation}
Here $\mu'$ is the \textit{transpose} of $\mu$. We give a more compact, diagrammatic description of the Tanisaki ideals in \Cref{sec:essential_generators}.

The cohomology rings $H^*(\overline{\cX_\mu})$ carry a non-obvious \cite{MR4029695} graded $\fS_n$-module structure, which is compatible with the natural action of $\fS_n$ on $\cR_\mu$. The graded \textit{Frobenius series} encoding the graded $\fS_n$-module decomposition of $\cR_\mu$ is the \textit{dual Hall--Littlewood symmetric function} up to a twist,
\begin{equation}\label{eq:Rmu_Qp}
  \GrFrob(\cR_\mu; q) = q^{b(\mu)} Q'_\mu(\mathbf{x}; q^{-1}) = \rev_q Q'_\mu(\mathbf{x}; q),
\end{equation}
where $b(\mu) \coloneqq \sum_i (i-1) \mu_i$ and $\rev_q F(q) \coloneqq q^{\deg F} F(q^{-1})$ is the \textit{$q$-reversal} operator. See \cite{MR1168926}. (The $\rev_q$ in \eqref{eq:Rmu_Qp} was inadvertently neglected in \cite[(7.1)]{MR3783430}.)

\subsection{Generalized coinvariant algebras}

In a different direction, Haglund--Rhoades--Shimozono \cite{MR3783430} introduced the \textit{generalized coinvariant algebras}
  \[ \cR_{n,k} \coloneqq \frac{\mathbb{Q}[x_1, \ldots, x_n]}{\langle x_1^k, \ldots, x_n^k, e_n(\underline{n}), e_{n-1}(\underline{n}), \ldots, e_{n-k+1}(\underline{n})\rangle} \]
while studying the \textit{Delta Conjecture} of Haglund--Remmel--Wilson \cite{MR3811519}, which we will discuss shortly. They gave the following compact description of the graded Frobenius series of $\cR_{n, k}$ \cite[Thm.~6.14]{MR3783430} generalizing \eqref{eq:Rmu_Qp} when $k=n$,
\begin{equation}\label{eq:Rnk_Qp}
  \GrFrob(\cR_{n, k}; q) = \rev_q \left[\sum_{\substack{\mu \vdash n \\ \ell(\mu) = k}} q^{\sum_{i=1}^k (i-1)(\mu_i - 1)} \binom{k}{m_1(\mu), \ldots, m_n(\mu)}_q Q_\mu'(\mathbf{x}; q)\right]
\end{equation}
where $\ell(\mu) \coloneqq \#\{j : \mu_j \neq 0\}$, $m_i(\mu) \coloneqq \#\{j : \mu_j = i\}$, and $\binom{k}{m_1, \ldots, m_n}_q$ is a \textit{$q$-multinomial coefficient}. 

Haglund--Rhoades--Shimozono ask whether a filtration of $\cR_{n, k}$ could be found to prove \eqref{eq:Rnk_Qp} directly using \eqref{eq:Rmu_Qp}, with successive quotients $\cR_\mu$ up to $q$-shifts \cite[Problem~7.1]{MR3783430}. A geometric description of $\cR_{n, k}$ was later given by Rhoades--Pawlowski \cite{MR4029695}, though an appropriate filtration has been elusive. Pursuing such a filtration has been the primary motivation of the present work.

\subsection{The Delta Conjecture and higher coinvariant algebras}

The Delta Conjecture of Haglund--Remmel--Wilson \cite{MR3811519} hypothesizes a certain symmetric function identity,
\begin{equation}\label{eq:Delta}
  \Delta_{e_{k-1}}'(e_n) = C_{n, k}(\mathbf{x}; q, t), \qquad 0 \leq k \leq n-1
\end{equation}
where $\Delta_f'$ is a certain modified Macdonald eigenoperator and $C_{n, k}(\mathbf{x}; q, t)$ is either of two explicit combinatorial expressions\footnote{The ``rise'' version has been independently proven by D'Adderio--Mellit \cite{MR4401822} and Blasiak--Haiman--Morse--Pun--Seelinger \cite{MR4553915}.}. See \cite[\S3]{MR3811519} for details.

The main result of Haglund--Rhoades--Shimozono is \cite[Thm.~6.11]{MR3783430},
\begin{equation}\label{eq:HRS.main}
  \GrFrob(\cR_{n, k}; q) = \rev_q \omega C_{n, k}(\mathbf{x}; q, 0),
\end{equation}
where $\omega \colon s_\lambda \mapsto s_{\lambda'}$ is the usual involution on symmetric functions. Representation-theoretically, $\omega$ corresponds to tensoring with the $\sgn$ representation. Consequently, $\cR_{n, k}$ provides a representation-theoretic model for the right-hand side of the $t=0$ specialization of the Delta Conjecture, up to a twist.

Zabrocki \cite{1902.08966} recently introduced the \textit{super-diagonal coinvariant algebra} $\cSDR_n$ and conjectured that it gives a representation-theoretic model for the left-hand side of the full Delta conjecture in the sense that
\begin{equation}\label{eq:Zabrocki}
  \GrFrob(\cSDR_n; q, t, z) = \sum_{k=0}^{n-1} z^k \Delta_{e_{n-k-1}}'(e_n).
\end{equation}
The $t=0$ specialization of Zabrocki's model differs from the generalized coinvariant algebras $\cR_{n, k}$ and is instead the \textit{super coinvariant algebra}
\begin{equation}
  \cSDR_n|_{t=0} = \cSR_n \coloneqq \frac{\mathbb{Q}[x_1, \ldots, x_n, \theta_1, \ldots, \theta_n]}{\cJ_n},
\end{equation}
where $\cJ_n$ is the ideal generated by the bi-homogeneous non-constant $\fS_n$-invariants. Here the $x_i$ commute, the $\theta_i$ \textit{anti-commute}, and $\fS_n$ acts simultaneously on $x$ and $\theta$ variables. That is, $x_i x_j = x_j x_i$, $x_i \theta_j = \theta_j x_i$, $\theta_i \theta_j = -\theta_j \theta_i$, $\sigma \cdot x_i = x_{\sigma(i)}$, and $\sigma \cdot \theta_i = \theta_{\sigma(i)}$. We may think of $\theta_{i_1} \cdots \theta_{i_k}$ as the differential $k$-form $\dif x_{i_1} \wedge \cdots \wedge \dif x_{i_k}$, and more generally $\mathbb{Q}[x_1, \ldots, x_n, \theta_1, \ldots, \theta_n]$ is the ring of differential forms on $V = \mathbb{Q}^n$ with polynomial coefficients.

The ideal $\cJ_n$ can be given a very explicit description. Let $\dif$ be the \textit{exterior derivative} defined by 
  \[ \dif f \coloneqq \sum_{i=1}^n \frac{\partial f}{\partial x_i} \dif x_i = \left(\sum_{i=1}^n \partial_{x_i} \theta_i\right) f. \]
Solomon showed \cite{MR0154929}
\begin{equation}
  \cJ_n = \langle e_r(\underline{n}), \dif e_r(\underline{n}) : r \in [n]\rangle.
\end{equation}

\begin{remark}
  After this work was submitted, Rhoades--Wilson \cite{RW23} proved the Hilbert series specialization of the $t=0$ case of Zabrocki's conjecture \eqref{eq:Zabrocki}. Consequently, the $p_1^n$ component of the formula \eqref{eq:SHI_Qp} below motivating this work has been entirely proven. Our results continue to provide additional evidence for \eqref{eq:SHI_Qp} and therefore for the full $t=0$ case of Zabrocki's conjecture.
\end{remark}

\subsection{Classical harmonics}

The coinvariant algebra $\cR_n$ has a distinguished set of coset representatives called the \textit{harmonics},
  \[ \cH_n \coloneqq \{f \in \mathbb{Q}[x_1, \ldots, x_n] : \partial_{e_r(\underline{n})} f = 0 \text{ for all }r \in [n]\}. \]
Here $\partial_g$ is the polynomial differential operator defined by replacing each $x_i$ with $\partial_{x_i}$. The natural projection $\cH_n \to \cR_n$ is an isomorphism of graded $\fS_n$-modules, so for many purposes we may replace $\cR_n$ with $\cH_n$. See \cite{sw2} for details. 

The alternating component of $\cH_n$,
  \[ \cH_n^{\sgn} \coloneqq \{f \in \cH_n : \forall \sigma \in \fS_n, \sigma \cdot f = \sgn(\sigma)f\}, \]
is spanned by the classical \textit{Vandermonde determinant},
  \[ \Delta_n \coloneqq \prod_{1 \leq i < j \leq n} (x_j - x_i). \]
Steinberg \cite[Thm.~1.3(c)]{MR167535} showed that
\begin{equation}\label{eq:Steinberg}
  \cH_n = \mathbb{Q}[\partial_{x_1}, \ldots, \partial_{x_n}] \Delta_n.
\end{equation}
Intuitively, we think of $\Delta_n$ as a ``tent pole'' which the remaining elements of $\cH_n$ ``hang off.''

\subsection{Harmonic differential forms}

Likewise, the super coinvariant algebras $\cSR_n$ may be replaced with the \textit{harmonic differential forms},
  \[ \cSH_n \coloneqq \{\omega \in \mathbb{Q}[x_1, \ldots, x_n, \theta_1, \ldots, \theta_n] : \partial_{e_r(\underline{n})} \omega = 0 = \partial_{\dif e_r(\underline{n})} \omega, r \in [n]\}. \]
Here $\partial_{\theta_i}$ is an \textit{interior product}. See \cite{sw2} for details. Let $\cSH_n^k$ denote the $k$-form component of $\cSH_n$.

In \cite{MR4258767}, Wallach and the author gave the following basis of size $2^{n-1}$ for the alternating component of $\cSH_n$,
\begin{align}
  \cSH_n^{\sgn}
    &\coloneqq \{\omega \in \cSH_n : \sigma \cdot \omega = \sgn(\sigma)\omega \text{ for all }\sigma \in \fS_n\} \nonumber \\
    &= \Span_{\mathbb{Q}}\{\dif_{i_1} \cdots \dif_{i_k} \Delta_n : 1 \leq i_1 < \cdots < i_k \leq n-1\}.
\end{align}
Here
  \[ \dif_i \coloneqq \sum_{j=1}^n \partial_{x_j}^i \theta_j \]
is a \textit{generalized exterior derivative} which lowers $x$-degree by $i$ and raises $\theta$-degree by $1$. For brevity, we write
  \[ \dif_I \coloneqq \dif_{i_1} \cdots \dif_{i_k} \]
where $I = \{i_1 < \cdots < i_k\} \subset [n-1]$. We sometimes abbreviate $\{i_1, \ldots, i_k\}$ as $i_1 \cdots i_k$.

In \cite{sw2}, Wallach and the author conjectured the following generalization of Steinberg's equation \eqref{eq:Steinberg},
\begin{equation}\label{eq:SHn_sgn}
  \cSH_n = \mathbb{Q}[\partial_{x_1}, \ldots, \partial_{x_n}] \{\dif_I \Delta_n : I \subset [n-1]\}.
\end{equation}
Among the evidence for \eqref{eq:SHn_sgn} provided in \cite{sw2}, we showed that the bi-graded support of $\cSH_n$ is precisely that predicted by \eqref{eq:SHn_sgn}, supporting the notion that the elements $\dif_I \Delta_n$ are the ``tent poles'' of $\cSH_n$. Rhoades--Wilson \cite{RW23} have since completely proven \eqref{eq:SHn_sgn}.

\begin{remark}
  Rhoades--Wilson \cite{MR4105531} have defined variations on the harmonics $\cSH_n$ by introducing ``superspace Vandermondes,'' which are alternants coming from particular terms in certain $\dif_I \Delta_n$'s. They construct modules by closing these superspace Vandermondes under partial derivatives which provably satisfy the appropriate analogue of \eqref{eq:SHI_Qp} below. It is an open problem to connect their modules to the harmonics $\cSH_n$.
\end{remark}

\subsection{The flip action}

Since $\cSH_n$ is closed under partial differentiation, we may consider it as a $\mathbb{Q}[x_1, \ldots, x_n]$-module under the \textit{flip action}
  \[ g \cdot \omega \coloneqq \partial_g \omega \]
for $\omega \in \cSH_n$. Since $e_d(\underline{n}) \cdot \omega = \partial_{e_d(\underline{n})} \omega = 0$ by definition, $\cSH_n$ is an $\cR_n$-module under the flip action. Note that the flip action \textit{lowers} $x$-degree.

Given $I \subset [n-1]$, define a component $\cSH_I$ of $\cSH_n$ from \eqref{eq:SHn_sgn} by
\begin{equation}
  \cSH_I \coloneqq \mathbb{Q}[\partial_{x_1}, \ldots, \partial_{x_n}] \dif_I \Delta_n.
\end{equation}
Let $\Ann \cSH_I$ be the annihilator of $\dif_I \Delta_n$ under the flip action, so that $\cSH_I \cong \mathbb{Q}[x_1, \ldots, x_n]/\!\Ann \cSH_I$ as $\mathbb{Q}[x_1, \ldots, x_n]$-modules.

Suppose for the sake of illustration that $\Ann \cSH_I = \cI_\mu$ is a Tanisaki ideal and $b(\mu) + i_1 + \cdots + i_k = \binom{n}{2}$. Since $\dif_{i_1} \cdots \dif_{i_k} \Delta_n$ transforms by $\sgn$ and has $x$-degree $\binom{n}{2} - i_1 - \cdots - i_k$, and since the flip action lowers $x$-degree, we have
\begin{align*}
  \GrFrob(\cSH_I; q)
    &= \omega q^{\binom{n}{2} - i_1 - \cdots - i_k} \GrFrob(\mathbb{Q}[x_1, \ldots, x_n]/\!\Ann \cSH_I; q^{-1}) \\
    &= \omega \rev_q \GrFrob(\mathbb{Q}[x_1, \ldots, x_n]/\cI_\mu; q) \\
    &= \omega Q_\mu'(\mathbf{x}; q).
\end{align*}
Consequently, using the super harmonics $\cSH_n$ and considering the flip action can entirely account for the twists in \eqref{eq:HRS.main}. We are thus led to the study of the $\mathbb{Q}[x_1, \ldots, x_n]$-module structure of $\cSH_n$.

\subsection{A potential filtration}\label{ssec:intro:filtration}

Combining \eqref{eq:SHn_sgn}, the $t=0$ case of Zabrocki's conjecture \eqref{eq:Zabrocki}, the $t=0$ case of the Delta Conjecture \eqref{eq:Delta}, and Haglund--Rhoades--Shimozono's $Q_\mu'$-expansion formula \eqref{eq:Rnk_Qp} gives
\begin{equation}\label{eq:SHI_Qp}
  \GrFrob\left(\sum_{\substack{I \subset [n-1]}} \cSH_I; q, z\right) = \sum_{\substack{\mu \vdash n}} z^{n-\ell(\mu)} q^{\sum_{i=1}^{\ell(\mu)} (i-1)(\mu_i - 1)} \binom{\ell(\mu)}{m_1(\mu), \ldots, m_n(\mu)}_q \omega Q_\mu'(\mathbf{x}; q).
\end{equation}

The left-hand side of \eqref{eq:SHI_Qp} is indexed by subsets of $[n-1]$. Expanding the multinomial coefficients, we may consider the right-hand side to be indexed by \textit{strong compositions} of $n$, namely sequences $\alpha = (\alpha_1, \ldots, \alpha_k)$ with $\alpha_i \geq 1$ and $\alpha_1 + \cdots + \alpha_k = n$, which are well-known to be in bijection with $2^{[n-1]}$. Combining all of these observations, we are led to the following question, which has motivated the present work. Here $\cI_\alpha \coloneqq \cI_\mu$ if $\mu$ is the weakly decreasing rearrangement of $\alpha$.

\begin{question}\label{que:main}
  Is there a total order $I_1 < I_2 < \cdots$ on $2^{[n-1]}$ and a bijection $\Phi_n$ from $2^{[n-1]}$ to the set of strong compositions $\alpha \vDash n$ for which the successive filtration quotients
    \[ \frac{\sum_{j \leq m} \cSH_{I_j}}{\sum_{j < m} \cSH_{I_j}} \]
  are annihilated precisely by the Tanisaki ideal $\cI_{\Phi_n(I_m)}$ acting as partial differential operators?
\end{question}

Additional motivation for considering \Cref{que:main} comes from a desire to find explicit bases for the super coinvariant algebras $\cSR_n$. Garsia--Procesi \cite{MR1168926} gave explicit monomial bases $\{x^\alpha\}$ for the Tanisaki ideals $\cI_\mu$. Hence given a total order and bijection satisfying \Cref{que:main}, we have an explicit basis for $\cSH_n$ of the form $\{\partial_{x^\alpha} \dif_I \Delta_n\}$.

We also have a purely enumerative consequence of \Cref{que:main}. In this situation,
\begin{equation}\label{eq:main_cor}
  \GrFrob\left(\sum_{\substack{I \subset [n-1]}} \cSH_I; q, z\right) = \sum_{I \subset [n-1]} z^{|I|} q^{\binom{n}{2} - \Sum(I) - b(\Phi_n(I))} \omega Q_{\Phi_n(I)}'(\mathbf{x}; q),
\end{equation}
where $\Sum(I) \coloneqq \sum_{i \in I} i$, $b(\alpha) \coloneqq b(\mu)$, and $Q_\alpha'(\mathbf{x}; q) \coloneqq Q_\mu'(\mathbf{x}; q)$ where $\mu$ is the weakly decreasing rearrangement of the strong composition $\alpha$. Define the \textit{coinversion number} of $\alpha \vDash n$ by
  \[ \coinv(\alpha) \coloneqq \#\{1 \leq i < j \leq \ell(\alpha) : \alpha_i < \alpha_j\}. \] 
Recall that
\begin{equation*}
  \sum q^{\coinv(\alpha)} = \binom{\ell(\mu)}{m_1(\mu), \ldots, m_n(\mu)}_q,
\end{equation*}
where the sum is over all rearrangements $\alpha$ of $\mu \vdash n$. Combining \eqref{eq:SHI_Qp} and \eqref{eq:main_cor} then gives
\begin{equation}\label{eq:Phi_condition}
  \sum_{\alpha \vDash n} z^{n-\ell(\alpha)} q^{2b(\alpha) - \binom{\ell(\alpha)}{2} + \coinv(\alpha)} = \sum_{I \subset [n-1]} z^{|I|} q^{\binom{n}{2} - \Sum(I)},
\end{equation}
where we have used the fact that $\omega Q_\mu'(\mathbf{x}; q)|_{\sgn} = q^{b(\mu)}$.

The classic ``stars and bars'' bijection from $2^{[n-1]}$ to $\{\alpha \vDash n\}$ does not satisfy \eqref{eq:Phi_condition}. In \Cref{sec:bijection}, we define a new bijection which does respect \eqref{eq:Phi_condition}. It is more convenient to describe the inverse map $\Psi_n$, which we do now. See \Cref{ex:Psi_n}.

\begin{definition}\label{def:Psi_n}
  Given a strong composition $\alpha$ of $n$, create a left-justified diagram of cells, where the $i$th row from the top has $\alpha_i$ cells. Let $m_i$ denote the number of cells in columns $1, 2, \ldots, i$. First fill the cells of the \textit{second} column from top to bottom with numbers $m_1, m_1-1, m_1-2, \ldots$, skipping missing cells in that column. Now delete the first column and any empty rows and repeat this procedure on the new second column using a maximum of $m_2$, and continue in this fashion. Afterwards, $\Psi_n(\alpha)$ is the set of numbers filling the columns $2, 3, \ldots$ of $\alpha$.
\end{definition}

\begin{example}\label{ex:Psi_n}
  When $\alpha = (1, 3, 2, 1, 3, 1) \vDash 11$, the procedure gives
  \begin{align*}
    \tableau{\ \\ \ & 5 & 9 \\ \ & 4 \\ \ \\ \ & 2 & 7 \\ \ }.
  \end{align*}
  Here $m_1 = 6$ and $m_2 = 9$. In the first phase, we fill the second column with numbers $6, 5, 4, 3, 2, 1$, skipping the missing cells $6, 3, 1$. In the second phase, we remove the first, fourth, and sixth rows and fill the remaining cells of the third column with $9, 8, 7$, skipping the missing cell $8$. In all, $\Psi_{11}(1, 3, 2, 1, 3, 1) = \{2, 4, 5, 7, 9\}$.
\end{example}

\begin{theorem}\label{thm:bijection}
  The bijection $\Psi_n \colon \{\alpha \vDash n\} \to 2^{[n-1]}$ satisfies
   \begin{align*}
     n-\ell(\alpha) &= |I| \\
     2b(\alpha) - \binom{\ell(\alpha)}{2} + \coinv(\alpha) &= \binom{n}{2} - \Sum(I)
   \end{align*}
   whenever $\Psi_n(\alpha) = I$. Consequently,
    \[  \sum_{\alpha \vDash n} z^{n-\ell(\alpha)} q^{2b(\alpha) - \binom{\ell(\alpha)}{2} + \coinv(\alpha)} = \sum_{I \subset [n-1]} z^{|I|} q^{\binom{n}{2} - \Sum(I)}. \]
\end{theorem}

\begin{remark}
  The condition \eqref{eq:Phi_condition} does not uniquely determine the bijection $\Phi_n$. For instance, one could replace $\coinv$ with $\inv$ or $\maj$ using a number of well-known bijections. As we show below, the bijection $\Phi_n$ from \Cref{thm:bijection} is sufficient to answer \Cref{que:main} for $1$-forms. However, computational evidence suggests a different order may be required in general. See \Cref{sec:further} for further discussion.
\end{remark}

\subsection{Tanisaki witness relations}

Given a total order and bijection satisfying \Cref{que:main}, for each generator $e_r(S)$ of the Tanisaki ideal $\cI_{\Phi_n(I_m)}$, we must have a relation of the form
\begin{equation}\label{eq:Tanisaki_witness}
  \partial_{e_r(S)} \dif_{I_m} \Delta_n = \sum_{j<m} \partial_{f_j} \dif_{I_j} \Delta_n \qquad\text{where}\qquad f_j \in \mathbb{Q}[x_1, \ldots, x_n],
\end{equation}
which we call a \textit{Tanisaki witness relation}. By homogeneity, we may restrict the terms in the Tanisaki witness relations to $k$-forms where $|I_j| = k$ is fixed.

\begin{example}
  When $n=3, k=1$, we have $\alpha \in \{(2, 1), (1, 2)\}$ with $I \in \{\{1\}, \{2\}\}$. The Tanisaki ideal $\cI_{(2, 1)}$ has the same generators as the classical coinvariant ideal $\cI_{(1^n)}$ together with $e_2(\underline{2})$ and its images under $\fS_3$. We find relations
  \begin{align*}
    \partial_{e_2(\underline{2})} \dif_{\{2\}} \Delta_3 &= 0 \\
    \partial_{e_2(\underline{2})} \dif_{\{1\}} \Delta_3 &= \partial_{e_1(\underline{2})} \dif_{\{2\}} \Delta_3.
  \end{align*}
  Hence $\cI_{(2, 1)}$ annihilates both $\cSH_{\{2\}} \subset \cSH_3^1$ and $(\cSH_{\{1\}}+\cSH_{\{2\}})/\cSH_{\{2\}}$, so the composition factors are both quotients of $\cI_{(2, 1)}$. By counting dimensions, there are no further relations, so the composition factors are precisely $\cI_{(2, 1)}$, answering \Cref{que:main} in the affirmative in this case using the order $\{2\} < \{1\}$.
\end{example}

\begin{example}
The relations between $\partial_{e_r(\underline{m})} \dif_I \Delta_n$ are generally quite complicated. For instance, at $n=7, k=2$, we have
\begin{align*}
  0 &= 5 \partial_{e_5(\underline{5})} \dif_{16} \Delta_7  - 4 \partial_{e_4(\underline{5})} \dif_{26} \Delta_7 + 3 \partial_{e_3(\underline{5})} \dif_{36} \Delta_7 - 2 \partial_{e_2(\underline{5})} \dif_{46} \Delta_7 +  \partial_{e_1(\underline{5})} \dif_{56} \Delta_7 \\
    &+ 3 \partial_{e_5(\underline{5})} \dif_{25} \Delta_7 - 2 \partial_{e_4(\underline{5})} \dif_{35} \Delta_7 + \partial_{e_3(\underline{5})} \dif_{45} \Delta_7 \\
    &+ \partial_{e_5(\underline{5})} \dif_{34} \Delta_7.
\end{align*}
and at $n=8, k=3$ we have
\begin{align*}
  0 &= 4 \partial_{e_6(\underline{6})} \dif_{356} \Delta_8
- 8 \partial_{e_5(\underline{6})} \dif_{357} \Delta_8
+ 4 \partial_{e_4(\underline{6})} \dif_{367} \Delta_8 \\
&- 3 \partial_{e_5(\underline{6})} \dif_{456} \Delta_8
+ 6 \partial_{e_4(\underline{6})} \dif_{457} \Delta_8
- 3 \partial_{e_3(\underline{6})} \dif_{467} \Delta_8.
\end{align*}
The first of these is explained by our results below, though the second is not.
\end{example}

\subsection{The generic Pieri rule}

All Tanisaki ideals $\cI_\mu$ with $\ell(\mu) = n-k$ for $\mu \neq (1^n)$ contain the generator $e_{n-k}(\underline{n-1})$. The following provides all necessary Tanisaki witness relations for this ``generic'' generator, and is one of our main results.

\begin{theorem}[``Generic Pieri Rule'']\label{thm:generic_Pieri}
  Suppose $I = \{i_1 < \cdots < i_k\} \subset [n-1]$.
  Then
    \[ \sum (-1)^d \partial_{e_{n-k-d}(\underline{n-1})} \dif_{j_1 \cdots j_k} \Delta_n = 0, \]
  where the sum is over all subsets $J = \{j_1 < \cdots < j_k\} \subset [n-1]$
  for which
    \[ 1 \leq i_1 \leq j_1 < i_2 \leq j_2 < i_3 \leq j_3 < \cdots < i_k \leq j_k < n, \]
  where
    \[ d \coloneqq (j_1 - i_1) + \cdots + (j_k - i_k). \]
\end{theorem}

\begin{remark}
  Our terminology in \Cref{thm:generic_Pieri} arises from the fact that the classical Pieri rule is a multiplicity-free expansion of the product of a Schur function by an elementary symmetric polynomial, together with the fact that the generator $e_{n-k}(\underline{n-1})$ is generic in the sense above.
\end{remark}

\begin{remark}\label{rem:generic}
  The generator $e_{n-k}(\underline{n-1})$, together with its images under $\fS_n$, is the only generator in $\cI_\mu$ for $\mu = (2^k, 1^{n-2k})$ when $n \geq 2k$, aside from the generators of $\cI_{(1^n)}$. The fraction of $\alpha \vDash n$ with $\ell(\alpha) = n-k$ where $\alpha$ is a rearrangement of $\mu = (2^k, 1^{n-2k})$ tends to $1$ for each fixed $k$ as $n \to \infty$. In this asymptotic sense, the Generic Pieri Rule gives ``almost all'' of the necessary Tanisaki witness relations.
\end{remark}

The Generic Pieri Rule answers the $1$-form case of \Cref{que:main} in the affirmative. More explicitly, we prove the following special case of \eqref{eq:main_cor}.

\begin{corollary}\label{cor:1_forms}
  The order $\{n-1\} < \{n-2\} < \cdots < \{1\}$ gives a filtration of $\cSH_n^1$ by $\cSH_{\{i\}}$'s where the composition factors are annihilated precisely by the Tanisaki ideal $\cI_{(2, 1^{n-2})}$. In particular,
  \begin{equation}\label{eq:cor:1_forms}
    \GrFrob\left(\cSH_n^1; q\right) = [n-1]_q \omega Q_{(2, 1^{n-2})}'(\mathbf{x}; q).
  \end{equation}
\end{corollary}

\subsection{Extreme hook relations}

In contrast to \Cref{thm:generic_Pieri}, which applies to any $I \subset [n-1]$, we also have Tanisaki witness relations corresponding to the least generic shapes in the following sense. Let $\overline{\alpha}$ be the result of removing the first column of $\alpha$ and removing empty rows, or equivalently subtracting $1$ from each entry and removing $0$'s. For $\alpha \vDash n$ with $\ell(\alpha) = n-k$, consider $\beta \coloneqq \overline{\alpha} \vDash k$. As noted above, for fixed $k$, the probability that $\beta = (1^k)$ tends to $1$ as $n \to \infty$. By contrast, the proportion of such $\alpha$ with $\beta = (k)$ is the smallest possible among all $\beta \vDash k$.

Slightly more generally, we consider $\alpha \vDash n$ with $\ell(\alpha) = n-k$ and $\overline{\alpha} = (s, 1^{k-s})$ for some $1 \leq s \leq k$. The Tanisaki ideal $\cI_\alpha$ is generated by
  \[ e_{n-k}(\underline{n-1}), e_{n-s}(\underline{n-2}), \ldots, e_{n-s}(\underline{n-s}) \]
together with their images under $\fS_n$ and the generators of the classical coinvariant ideal $\cI_{(1^n)}$. The following result gives Tanisaki witness relations for each of these generators.

\begin{theorem}\label{thm:hook_witness_extreme}
  Suppose $I = \{i_1 < \cdots < i_k\} \subset [n-1]$ is such that for some $1 \leq s \leq k$ we have
  \begin{align*}
    i_1, \ldots, i_{k-s+1} &\leq n-k \\
    i_{k-s+2} &= n-s+1 \\
    i_{k-s+3} &= n-s+2 \\
    &\ \,\vdots \\
    i_k &= n-1.
  \end{align*}
  Pick $0 \leq u \leq s$. Then
  \begin{equation}\label{eq:hook_witness_extreme.1}
    \sum (-1)^d \Delta_s(j_{k-s+1}, \ldots, j_k) \binom{d+u}{u}
        \partial_{e_{n-s-d}(\underline{n-s+u})} \dif_J \Delta_n = 0,
  \end{equation}
  where the sum is over all subsets $J = \{j_1 < \cdots < j_k\} \subset [n-1]$
  for which
  \begin{align*}
    j_1 = i_1, \ldots, j_{k-s} &= i_{k-s} \\
    d \coloneqq (j_{k-s+1} - i_{k-s+1}) + &\cdots + (j_k - i_k) \geq 0.
  \end{align*}
\end{theorem}

\begin{remark}\label{rem:extreme_hook}
  The condition on $I$ in \Cref{thm:hook_witness_extreme} is equivalent to $\Phi_n(I) = \alpha$ where $\alpha \vDash n$, $\ell(\alpha) = n-k$, and $\overline{\alpha} = (s, 1^{k-s})$ for some $1 \leq s \leq k$.
\end{remark}

\subsection{Paper organization}

The rest of the paper is organized as follows. In \Cref{sec:essential_generators}, we describe a set of ``essential'' Tanisaki ideal generators. In \Cref{sec:bijection}, we give the inverse to the bijection $\Psi_n$ from \Cref{def:Psi_n} and prove \Cref{thm:bijection}. In \Cref{sec:staircase}, we introduce a combinatorial model for the terms in our main identities. In \Cref{sec:Pieri}, we prove the Generic Pieri Rule, \Cref{thm:generic_Pieri}, and the $1$-form result, \Cref{cor:1_forms}. In \Cref{sec:actions}, we introduce some symmetric group actions and give a shifted Vandermonde identity, \Cref{cor:shifted_Vandermonde}. In \Cref{sec:hook}, we use the results of the previous sections to prove our second family of Tanisaki witness relations, \Cref{thm:hook_witness_extreme}. In \Cref{sec:further}, we discuss further directions.

\section{Essential Tanisaki generators}\label{sec:essential_generators}

We now describe a small subset of the Tanisaki ideal generators which in fact suffice to generate $\cI_\mu$. See \Cref{ex:essential_Tanisaki} for a simple graphical interpretation of this set of ``essential'' generators.

\begin{lemma}
  Given $\mu \vdash n$, compute $\underline{d}_0, \ldots, \underline{d}_{\mu_1-1}$ iteratively by $\underline{d}_0 \coloneqq 1$ and
    \[ \underline{d}_i \coloneqq \underline{d}_{i-1} + (\mu_i' - 1). \]
  Then
  \begin{equation}\label{eq:cI_mu}
    \cI_\mu = \fS_n \cdot \langle e_{\underline{d}_1}(\underline{n-1}), e_{\underline{d}_2}(\underline{n-2}), \ldots, e_{\underline{d}_{\mu_1-1}}(\underline{n-\mu_1+1}), e_1(\underline{n}), \ldots, e_n(\underline{n})\rangle.
  \end{equation}

  \begin{proof}
    First recall from \eqref{eq:Tanisaki_generators} that the Tanisaki ideal associated to $\mu \vdash n$ is by definition
      \[ \cI_\mu \coloneqq \langle T_\mu\rangle \]
    where
      \[ T_\mu \coloneqq \{e_r(S) : |S| - d_{|S|}(\mu) < r \leq |S|, S \subset [n]\} \]
    with $d_k(\mu) \coloneqq \mu_n' + \mu_{n-1}' + \cdots + \mu_{n-k+1}'$. Here $\mu'$ is padded with $0$'s if necessary so that it has $n$ entries.

    We have $n-d_n(\mu) = 0$, so $e_r(\underline{n}) \in T_\mu$ for $1 \leq r \leq n$. We similarly have $e_{\underline{d}_i}(\underline{n-i}) \in T_\mu$ for $1 \leq i \leq \mu_1 - 1$ if
      \[ (n-i) - d_{n-i} < \underline{d}_i, \]
    or equivalently if
      \[ (n-i) - \mu_n' - \cdots - \mu_{i+1}' \leq \mu_1' + \cdots + \mu_i' - i. \]
    Equality holds in this last expression, so in fact $(n-i) - d_{n-i} = \underline{d}_i - 1$. Write $\cI_\mu'$ for the right-hand side of \eqref{eq:cI_mu}. We have just shown that $\cI_\mu'$ is contained in $\cI_\mu$.
    
    Conversely, we show $e_r(S) \in \cI_\mu'$ for $|S| - d_{|S|}(\mu) < r \leq |S|$, $S \subset [n]$ by downward induction on $|S|$. By $\fS_n$-symmetry, we may suppose $S = \{1, 2, \ldots, n-i\}$. In the base case $i=0$, $e_r(\underline{n}) \in \cI_\mu'$. Next suppose $0<i<\mu_1$. We further induct on $r$. In the base case, $r = (n-i) - d_{n-i}(\mu) + 1 = \underline{d}_i$ and $e_r(\underline{n-i}) \in \cI_\mu'$ by assumption. For $r > \underline{d}_i$, we have the simple identity
      \[ e_r(\underline{n-i+1}) = e_r(\underline{n-i}) + x_{n-i+1} e_{r-1}(\underline{n-i}). \]
    By induction on $r$, $e_{r-1}(\underline{n-i}) \in \cI_\mu'$. On the other hand, $r > \underline{d}_i \geq \underline{d}_{i-1}$, so $e_r(\underline{n-(i-1)}) \in \cI_\mu'$ by induction on $i$. Hence $e_r(\underline{n-i}) \in \cI_\mu'$, completing the induction on $r$, and hence on $i$. Finally, if $i \geq \mu_1$, we have $d_{n-i} = 0$, so no such $r$ exists, completing the proof.
  \end{proof}
\end{lemma}

\begin{example}\label{ex:essential_Tanisaki}
  Let $\mu = (5, 3, 1, 1, 1)$. After drawing the diagram of $\mu$, compute the sequence $\underline{d}$ by writing $1$ above the first column, adding one less than the length of the first column and writing the result above the second column, etc. Here we have
  \begin{align*}
    \underline{d} &= 1,\, 5,\, 6,\, 7,\, 7 \\
    &\quad\tableau{\ & \ & \ & \ & \ \\ \ & \ & \ \\ \ \\ \ \\ \ }
  \end{align*}
  so that
    \[ \cI_{(5, 3, 1, 1, 1)} = \fS_n \cdot \langle e_5(\underline{n-1}), e_6(\underline{n-2}), e_7(\underline{n-3}), e_7(\underline{n-4}), e_1(\underline{n}), \ldots, e_n(\underline{n})\rangle \]
  where $n=11$.
\end{example}

\begin{example}
  Suppose $\mu \vdash n$ has $\ell(\mu) = n-k$ and $\mu \neq (1^n)$. Then $\underline{d}_1 = n-k$, so $e_{n-k}(\underline{n-1}) \in \cI_\mu$, which is the ``generic'' generator involved in the Generic Pieri Rule, \Cref{thm:generic_Pieri}. Moreover, if $\overline{\mu} = (1^k)$, so $\mu = (2^k, 1^{n-2k})$, this is the only generator up to the $\fS_n$-action aside from the generators of the classical coinvariant ideal $\cI_{(1^n)}$.
\end{example}

\begin{example}
  Suppose $\mu \vdash n$ has $\ell(\mu) = n-k$ and $\overline{\mu} = (s, 1^{k-s})$ for $1 \leq s \leq k$. Then $\underline{d} = 1, n-k, n-s, \ldots, n-s$ and the essential generators of $\cI_\mu$ are
    \[ e_{n-k}(\underline{n-1}), e_{n-s}(\underline{n-2}), \ldots, e_{n-s}(\underline{n-s}) \]
\end{example}

\section{Subset to composition bijection}\label{sec:bijection}

We now describe the inverse $\Phi_n \colon 2^{[n-1]} \to \{\alpha \vDash n\}$ to the map $\Psi_n$ from \Cref{ssec:intro:filtration} described in \Cref{def:Psi_n}. Along the way, we prove the statistic preservation result for the maps $\Phi_n$ and $\Psi_n$, \Cref{thm:bijection}. This section may be read independently of the others.

We begin by considering a step of a recursive decomposition on strong compositions. We also define a notion of ``degree'' inspired by \eqref{eq:Phi_condition} and describe the effect of this recursive step on the degree.

\begin{definition}
  Let $\alpha \vDash n$ be a strong composition of $n$ of length $\ell(\alpha)$. Set
    \[ \coinv(\alpha) \coloneqq \#\{1 \leq i < j \leq \ell(\alpha) : \alpha_i < \alpha_j\}. \]
  Let $\mu(\alpha)$ denote the partition of $n$ obtained by rearranging $\alpha$ in weakly decreasing order. Set
    \[ \deg(\alpha) \coloneqq \coinv(\alpha) + \sum_{i=1}^{\ell(\alpha)} (i-1)(2\mu(\alpha)_i - 1) = \coinv(\alpha) + 2b(\alpha) - \binom{\ell(\alpha)}{2}. \]
  Finally, let $\overline{\alpha}$ be the strong composition obtained by removing $1$ from every row of $\alpha$ and deleting empty rows.
\end{definition}

\begin{example}
  When $\alpha = (1, 3, 2, 1, 3, 1)$, we have $\mu(\alpha) = (3, 3, 2, 1, 1, 1)$, $\overline{\alpha} = (2, 1, 2)$ so
    \[ \alpha = \tableau{\ & \\ \ & \ & \ \\ \ & \ \\ \ \\ \ & \ & \ \\ \ }
        \qquad \mu(\alpha) = \tableau{\ & \ & \ \\ \ & \ & \ \\ \ & \ \\ \ \\ \ \\ \ }
        \qquad \overline{\alpha} = \tableau{\ & \ \\ \ \\ \ & \ } \]
  and
  \begin{align*}
    \coinv(\alpha) &= 3 + 0 + 1 + 1 + 0 + 0 = 5 \\
    b(\alpha) &= 0 \cdot 3 + 1 \cdot 3 + 2 \cdot 2 + 3 \cdot 1 + 4 \cdot 1 + 5 \cdot 1 = 19 \\
    \deg(\alpha) &= 5 + 2 \cdot 19 - 15 = 28 \\
    \deg(\overline{\alpha}) &= 1 + 2 \cdot 4 - 3 = 6.
  \end{align*}
\end{example}

\begin{lemma}\label{lem:alphabar}
  Let $\alpha \vDash n$. Suppose $\ell(\alpha) = r$ and $\ell(\overline{\alpha}) = s$. Then
    \[ \deg(\alpha) = \deg(\overline{\alpha}) + \binom{r}{2} - s + j_1 + \cdots + j_s \]
  where $\{1 \leq j \leq r : \alpha_j > 1\} = \{j_1, \ldots, j_s\}$.
  
  \begin{proof}
    By considering coinversions of $\alpha$ starting from a row of length $1$ separately, it is easy to see that
      \[ \coinv(\alpha) = \coinv(\overline{\alpha}) + (j_1 - 1) + \cdots + (j_s - s). \]
    On the other hand, we have
    \begin{align*}
      &\sum_{i=1}^r (i-1)(2\mu(\alpha)_i - 1)
             - \sum_{i=1}^s (i-1)(2\mu(\overline{\alpha})_i - 1) \\
        &= \sum_{i=s+1}^r (i-1)(2 \cdot 1 - 1)
             + \sum_{i=1}^s (i-1) 2(\mu(\alpha)_i - \mu(\overline{\alpha})_i) \\
        &= \sum_{i=s+1}^r (i-1) + \sum_{i=1}^s 2(i-1) \\
        &= \binom{r}{2} + \sum_{i=1}^s (i-1).
    \end{align*}
    The result follows by combining these observations.
  \end{proof}
\end{lemma}


We likewise consider a step of a recursive decomposition on subsets of $[n-1]$. We again define a notion of ``degree'' inspired by \eqref{eq:Phi_condition} and describe the effect of this recursive step on the degree. Finally we restate and prove \Cref{thm:bijection}.

\begin{definition}
  Fix $n$. Let $I \subset [n-1]$. Define
    \[ \deg(I) \coloneqq \binom{n}{2} - \sum_{i \in I} i. \]
  Suppose $I = \{i_1 < \cdots < i_k\}$. Let $\overline{I} \subset [k-1]$ be defined as follows. There is some unique $1 \leq s \leq k$ such that
  \begin{equation}\label{eq:is_nk}
    1 \leq i_1 < \cdots < i_s \leq n-k < i_{s+1} < \cdots < i_k \leq n-1.
  \end{equation}
  Set
    \[ \overline{I} \coloneqq \{i_1', \ldots, i_{k-s}'\} \subset [k-1] \qquad \text{ where }
       \qquad i_j' \coloneqq i_{s+j} - n + k. \]
\end{definition}

\begin{example}
  Let $n=11$ and $I = \{2, 4, 5, 7, 9\} \subset [10]$. Here $k=5$ and
    \[ i_1 < i_2 < i_3 \leq n-k = 6 < i_4 < i_5, \]
  so $s=3$ and $\overline{I} = \{i_4 - 6, i_5 - 6\} = \{7 - 6, 9 - 6\} = \{1, 3\} \subset [4]$.
  We see
  \begin{align*}
    \deg(I)
      &= \binom{11}{2} - (2+4+5+7+9) = 28 \\
    \deg(\overline{I})
      &= \binom{5}{2} - (1+3) = 6.
  \end{align*}
\end{example}

\begin{lemma}\label{lem:Ibar}
  Let $I = \{i_1 < \cdots < i_k\} \subset [n-1]$ and
  $\overline{I} = \{i_1', \ldots, i_{k-s}'\} \subset [k-1]$ as above. Then
    \[ \deg(I) = \deg(\overline{I}) + \binom{n-k}{2} + s(n-k) - \sum_{j=1}^s i_j. \]
  
  \begin{proof}
    We compute
    \begin{align*}
      \deg(I) - \deg(\overline{I})
        &= \binom{n}{2} - \binom{k}{2} + \sum_{j=1}^{k-s} i_j' - \sum_{j=1}^k i_j \\
        &= \binom{n}{2} - \binom{k}{2} + \sum_{j=1}^{k-s} (i_{s+j} - n + k) - \sum_{j=1}^k i_j \\
        &= \binom{n}{2} - \binom{k}{2} - (n-k)(k-s) + \sum_{j=s+1}^k i_j - \sum_{j=1}^k i_j \\
        &= \binom{n-k}{2} + s(n-k) - \sum_{j=1}^s i_j.
    \end{align*}
  \end{proof}
\end{lemma}

\begin{definition}
  We recursively define a bijection $\Phi_n$ from subsets $I$ of $[n-1]$ to strong compositions $\alpha$ of $n$ as follows. Take $|I| = k$. We will ensure $\ell(\Phi_n(I)) = n-k$. For $k=0$, set $\Phi_n(\varnothing) = (1^n)$. For $k>0$, we have $\overline{I} \subset [k-1]$ and $s$ satisfying \eqref{eq:is_nk}. Let $\beta = \Phi_{k}(\overline{I})$, so $\ell(\beta) = n-k-|\overline{I}| = k-(k-s) = s$. Construct $\alpha$ from $\beta$ by requiring $\overline{\alpha} = \beta$ and
    \[ \{1 \leq j \leq n-k : \alpha_j > 1\} = \{j_1 < \cdots < j_s\} \]
  where
  \begin{align*}
    j_1 = n-k&+1 - i_s \\
    &\vdots \\
    j_s =n-k&+1 - i_1.
  \end{align*}
\end{definition}

\begin{example}
  Consider
  \begin{alignat*}{2}
    &I = \{2, 4, 5 \mid 7, 9\} \subset [11-1] \qquad
      &&\alpha = (1, 3, 2, 1, 3, 1) \\
    &\overline{I} = \{1, 3 \mid\ \} \subset [5-1]
      &&\overline{\alpha} = (2, 1, 2) \\
    &\overline{\overline{I}} = \varnothing \subset [2-1]
      &&\overline{\overline{\alpha}} = (1, 1)
  \end{alignat*}
  Here $\mid$ indicates the two halves of the decompositions from \eqref{eq:is_nk}. The corresponding diagrams using the inverse map $\Psi_n$ from \Cref{def:Psi_n} are
  \begin{align*}
    \tableau{\ \\
                   \ } \qquad
    \tableau{\ & 3\\
                   \ \\
                   \ & 1} \qquad
    \tableau{\ \\
                   \ & 5 & 9 \\
                   \ & 4 \\
                   \ \\
                   \ & 2 & 7 \\
                   \ }
  \end{align*}
  The elements left of $\mid$ indicate where to attach elements of $\overline{\alpha}$ to $(1^{n-k})$ to form $\alpha$, from right to left. We have $\Phi_2(\overline{\overline{I}}) = \overline{\overline{\alpha}}$, $\Phi_5(\overline{I}) = \overline{\alpha}$, and $\Phi_{11}(I) = \alpha$.
\end{example}

\begin{reptheorem}{thm:bijection}
  The bijection $\Psi_n \colon \{\alpha \vDash n\} \to 2^{[n-1]}$ satisfies
   \begin{align*}
     n-\ell(\alpha) &= |I| \\
     2b(\alpha) - \binom{\ell(\alpha)}{2} + \coinv(\alpha) &= \binom{n}{2} - \Sum(I)
   \end{align*}
   whenever $\Psi_n(\alpha) = I$. Consequently,
    \[  \sum_{\alpha \vDash n} z^{n-\ell(\alpha)} q^{2b(\alpha) - \binom{\ell(\alpha)}{2} + \coinv(\alpha)} = \sum_{I \subset [n-1]} z^{|I|} q^{\binom{n}{2} - \Sum(I)}. \]
   
  \begin{proof}
    We've ensured $n-\ell(\alpha) = |I|$. The second condition is equivalent to
      \[ \deg(\Phi_n(I)) = \deg(I). \]
    In the base case,
      \[ \deg(\Phi_n(\varnothing)) = \deg((1^n)) = \sum_{i=1}^n (i-1) = \binom{n}{2}
          = \deg(\varnothing). \]
    Inductively, we may suppose that $\deg(\overline{I}) = \deg(\overline{\alpha})$.
    By \Cref{lem:alphabar} and \Cref{lem:Ibar}, where $r = n-k$,
    \begin{align*}
      \deg(I)
        &= \deg(\overline{I}) + \binom{n-k}{2} + s(n-k) - (i_1 + \cdots + i_s) \\
        &= \deg(\overline{I}) + \binom{r}{2} + s(n-k) - (s(n-k+1) - j_s - \cdots - j_1) \\
        &= \deg(\overline{\alpha}) + \binom{r}{2} - s + j_1 + \cdots + j_s \\
        &= \deg(\alpha),
    \end{align*}
    which completes the proof.
  \end{proof}
\end{reptheorem}

We also note that, from this recursive description, it is easy to see that $\Phi_n$ and $\Psi_n$ are in fact inverses, hence bijections.

\section{Marked staircase diagrams}\label{sec:staircase}

We now introduce a combinatorial model for the terms in $\partial_{e_r(\underline{m})} \dif_I \Delta_n$ using decorated diagrams. We will use relations between these diagrams to build sign-reversing involutions in the subsequent sections.

\subsection{Staircases}

Let $\Delta_n \coloneqq \prod_{1 \leq i < j \leq n} (x_j - x_i)$ denote the Vandermonde determinant in $n$ variables. We have
\begin{equation}\label{eq:Vandermonde_expanded}
  \Delta_n = \sum_{\sigma \in S_n} (-1)^{\sgn(\sigma)} x_1^{\sigma(1) - 1} \cdots x_n^{\sigma(n) - 1}.
\end{equation}
We model the monomials appearing in \eqref{eq:Vandermonde_expanded} as follows.

\begin{definition}
  An \textit{$n$-staircase} is a bottom-justified arrangement of $n$ columns of cells with heights $0, 1, \ldots, n-1$, each used exactly once. The \textit{sign} of an $n$-staircase with column heights $h_1, \ldots, h_n$ in order from left to right is $(-1)^c$ where
    \[ c = \#\{i < j : h_i > h_j\}. \]
  Equivalently, the sign is $\sgn \Delta_n(h_1, \ldots, h_n)$, where we have used the signum function. The \textit{monomial weight} of such an $n$-staircase is $x_1^{h_1} \cdots x_n^{h_n}$, and the \textit{weight} is $(-1)^c x_1^{h_1} \cdots x_n^{h_n}$.
\end{definition}

\begin{example}
  The $6$-staircase with heights $h_1=1, h_2=5, h_3=3, h_4=0, h_5=2, h_6=4$ is
  
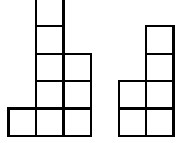
\[ \begin{ytableau}
        \none & \ \\
        \none & \ & \none & \none & \none & \ \\
        \none & \ & \ & \none & \none & \ \\
        \none & \ & \ & \none & \ & \  \\
        \ & \ & \ & \none & \ & \ 
      \end{ytableau} \]
  and has weight $(-1)^7 x_1 x_2^5 x_3^3 x_5^2 x_6^4$.
\end{example}

By \eqref{eq:Vandermonde_expanded}, $\Delta_n$ is the weight generating function of the $n$-staircases.

\subsection{Marked staircases}

Monomials in $\partial_{e_r(\underline{m})} \dif_I \Delta_n$ arise from applying some sequence of operators $\partial_{x_j}^{i_j} \theta_j$ to a monomial from \eqref{eq:Vandermonde_expanded}, followed by $\partial_{x_J}$ for some $J \subset [m]$. We model these terms diagrammatically as follows. See \Cref{ex:marked_staircase}.

\begin{definition}\label{def:marked_staircases}
  A \textit{marked staircase} is an $n$-staircase where some of the boxes have been filled with $\times$'s or $\circ$'s subject to the following constraints:
  \begin{enumerate}
    \item Any $\times$'s are top-justified in their column.
    \item Any $\circ$'s are top-justified in their column below any $\times$'s.
    \item A column may have at most one $\circ$.
    \item The last $i \geq 0$ columns are colored grey and are forbidden from containing $\circ$'s. They may still contain $\times$'s.
  \end{enumerate}
  Furthermore, the \textit{weight} of a marked staircase is the product of the following three terms.
  \begin{itemize}
    \item The \textit{monomial weight} of a marked staircase is $x_1^{g_1} \cdots x_n^{g_n} \theta_{c_1} \cdots \theta_{c_k}$ where $g_\ell$ denotes the number of \textit{unmarked} boxes in column $\ell$ and $\{c_1 < \cdots < c_k\}$ is the set of indexes of columns which contain $\times$'s.
    \item The \textit{sign} of a marked staircase is $(-1)^c \sgn \Delta_k(j_1, \ldots, j_k)$ where $(-1)^c$ is the sign of the underlying $n$-staircase and $j_\ell$ is the number of $\times$'s in column $c_\ell$. Note that this is zero if and only if $j_1, \ldots, j_k$ are not all distinct.
    \item The \textit{order} of a marked staircase is the product of the heights at which the $\times$'s and $\circ$'s appear.
  \end{itemize}
\end{definition}

\begin{example}\label{ex:marked_staircase}
  The marked $6$-staircase
  
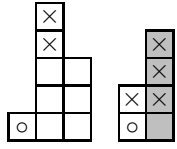
\[ \begin{ytableau}
        \none & \times \\
        \none & \times & \none & \none & \none & *(lightgray)\times \\
        \none & \ & \ & \none & \none & *(lightgray)\times \\
        \none & \ & \ & \none & \times & *(lightgray)\times \\
        \circ & \ & \ & \none & \circ & *(lightgray)
      \end{ytableau} \]
  has monomial weight $x_2^3 x_3^3 x_6 \theta_2 \theta_5 \theta_6$, sign $(-1)^7 \sgn \Delta_3(2, 1, 3) = 1$, and order $1 \cdot (5 \cdot 4) \cdot (2 \cdot 1) \cdot (4 \cdot 3 \cdot 2) = 960$. The weight is thus $960 x_2^3 x_3^3 x_6 \theta_2 \theta_5 \theta_6$, which represents a term in $\partial_{e_2(\underline{5})} \dif_{123} \Delta_6$.
\end{example}

\begin{lemma}\label{lem:marked_staircases_gf}
  Suppose $I = \{i_1 < \cdots < i_k\} \subset [n-1]$. Then
    \[ \partial_{e_r(\underline{n-m})} \dif_I \Delta_n \]
  is the weight generating function for marked $n$-staircases with $r$ $\circ$'s, $\times$'s of lengths $i_1, \ldots, i_k$, and the last $m$ columns grey.
  
  \begin{proof}
    Applying $\dif_{i_\ell}$ to $\Delta_n$ is essentially the same as picking a marked $n$-staircase and picking a column to add $i_\ell$ $\times$'s to, ignoring scalars and the $\theta$-part for the moment. Analogously, applying $\partial_{e_r(\underline{n-m})}$ is the same as picking $r$ of the first $n-m$ columns to add $\circ$'s to, namely the non-grey columns. The scalars arising from applying these $x$-derivatives are precisely the product of the heights of the marks involved, which is the order. The $x$-part of the monomial weight is thus correct.
    
    For the $\theta$-part, suppose $c_\ell'$ is the index of the column with $i_\ell$ $\times$'s. We are hence tracking the term $\partial_{x_{c_1'}}^{i_1} \theta_{c_1'} \cdots \partial_{x_{c_k'}}^{i_k} \theta_{c_k'}$ in $\dif_{i_1} \cdots \dif_{i_k}$, so the required $\theta$-part is $\theta_{c_1'} \cdots \theta_{c_k'}$. Let $c_1 < \cdots < c_k$ be the increasing rearrangement of $c_1', \ldots, c_k'$ and say that column $c_\ell$ has $j_\ell$ $\times$'s. Let $c_\ell = c_{\sigma(\ell)}'$ for some $\sigma \in \fS_k$, so that $j_\ell = i_{\sigma(\ell)}$. We have
    \begin{align*}
      \theta_{c_1'} \cdots \theta_{c_k'}
       &= \sgn(\sigma) \theta_{c_1} \cdots \theta_{c_k}.
    \end{align*}
    Since $i_1 < \cdots < i_k$ has the same relative order as $1 < \cdots < k$,
    \begin{align*}
      \sgn(\sigma)
        &= \sgn \Delta_n(\sigma(1), \ldots, \sigma(k)) \\
        &= \sgn \Delta_n(i_{\sigma(1)}, \ldots, i_{\sigma(k)}) \\
        &= \sgn \Delta_n(j_1, \ldots, j_k).
    \end{align*}
  \end{proof}
\end{lemma}

\begin{remark}
  \Cref{lem:marked_staircases_gf} remains valid if we use a multiset $\{\!\{i_1 \leq \cdots \leq i_k\}\!\} \subset [n-1]$, since if the indexes are not all distinct, $\dif_I = 0$ and the weights are zero. We will use such degenerate terms in a later argument.
\end{remark}

\subsection{Marked staircase relations}

The following operations preserve or negate the monomial weight of a marked staircase.  We provide examples of each operation, where the altered portions have been highlighted.

\begin{lemma}\label{lem:relations}
  \ 
  \begin{enumerate}[A.]
    \item Taking a non-grey column with at least two $\times$'s and without an $\circ$ and replacing the bottommost $\times$ with an $\circ$ toggles the parity of the number of $\circ$'s and preserves weight if it is non-zero.
      \begin{align*}
        \begin{ytableau}
          \none & \none & \times & \none \\
          \none & \none & \times & *(lightgray) \times \\
          \circ & \none & *(yellow) \times & *(lightgray)
        \end{ytableau}
        \quad=\quad
        \begin{ytableau}
          \none & \none & \times & \none \\
          \none & \none & \times & *(lightgray) \times \\
          \circ & \none & *(yellow) \circ & *(lightgray)
        \end{ytableau}
      \end{align*}
    \item Taking non-grey columns of height $v$ and $v-1$ where the column of height $v$ has an $\circ$ and the column of height $v-1$ does not have an $\circ$ and swapping the columns and the $\circ$ while preserving the number of $\times$'s in each original column negates the weight and preserves the number of $\circ$'s.
      \begin{align*}
        \begin{ytableau}
          \none & \none & \none & *(lightgray) \times \\
          \none & \none & *(lime) \times & *(lightgray) \times \\
          \ & \none & *(yellow) \circ & *(lightgray) \times
        \end{ytableau}
        \quad=\quad
        -1 \cdot \begin{ytableau}
          \none & \none & \none & *(lightgray) \times \\
          *(yellow) \circ & \none & \none & *(lightgray) \times \\
          \ & \none & *(lime) \times & *(lightgray) \times
        \end{ytableau}
      \end{align*}
    \item If $j$ $\times$'s appear in a non-grey column with no $\circ$ and $j-1$ $\times$'s appear in a column with an $\circ$, swapping the final $\times$ and $\circ$ negates the weight and preserves the number of $\circ$'s, assuming $j \geq 2$.
      \begin{align*}
        \begin{ytableau}
          \none & \none & \times & \none \\
          \times & \none & *(yellow)\times & \none \\
          *(lime)\circ & \none & \ & *(lightgray)\ 
        \end{ytableau}
        \quad=\quad
        -1 \cdot
        \begin{ytableau}
          \none & \none & \times & \none \\
          \times & \none & *(lime)\circ & \none \\
          *(yellow)\times & \none & \ & *(lightgray)\ 
        \end{ytableau}
      \end{align*}
    \item Given two columns with blocks of $\times$'s at the same height, we may move the stack of $\times$'s above the common height from one column to the other. The number of $\circ$'s is preserved and the weight is either preserved up to a sign or is zero.
      \begin{align*}
        \begin{ytableau}
          \none & \none & \none & *(cyan)\times & \none & \none \\
          \ & \none & \none & *(cyan)\times & \none & \none \\
          \ & *(pink)\times & \none & *(pink)\times & \none & \none \\
          \ & *(yellow)\ & \none & *(lime)\times & \none & *(lightgray)\ \\
          \ & *(yellow)\ & \none & *(lime)\ & *(lightgray)\ & *(lightgray)\ 
        \end{ytableau}
        \quad=\quad
        \pm 1 \cdot
        \begin{ytableau}
          \none & *(cyan)\times & \none & \none & \none & \none \\
          \ & *(cyan)\times & \none & \none & \none & \none \\
          \ & *(pink)\times & \none & *(pink)\times & \none & \none \\
          \ & *(yellow)\ & \none & *(lime)\times & \none & *(lightgray)\ \\
          \ & *(yellow)\ & \none & *(lime)\ & *(lightgray)\ & *(lightgray)\ 
        \end{ytableau}
      \end{align*}
  \end{enumerate}
\end{lemma}

\noindent In the particular example for relation (D), the parity of $c$ changes, which cancels with the sign change from the $\theta$-part, so the weight is in fact preserved.

\section{Generic Pieri Rule proof}\label{sec:Pieri}

We now turn to the proof of our first, larger family of Tanisaki witness relations, \Cref{thm:generic_Pieri}. Our overall strategy will be to collect together certain types of marked staircases and cancel them amongst themselves using relations (A)-(C) from \Cref{lem:relations}. Before proving \Cref{thm:generic_Pieri}, we introduce some notation used in the proof.

\begin{notation}\label{not:generic_Pieri}
  For $J = \{j_1 < \cdots < j_k\} \subset [n-1]$, let $\{j_1, \ldots, j_k\}^r$ denote the weight generating function of the marked $n$-staircases with $\times$'s of lengths $j_1 < \cdots < j_k$, $r$ $\circ$'s, and the final column greyed out. By \Cref{lem:marked_staircases_gf},
    \[ \{j_1, \ldots, j_k\}^r = \partial_{e_r(\underline{n-1})} \dif_J \Delta_n. \]
  Additionally, we decorate $j_1, \ldots, j_k$ to indicate the weight generating function of such staircases subject to the following mutually exclusive and exhaustive constraints:
  \begin{enumerate}[(i)]
    \item $j^\flat$ means the column with $j$ $\times$'s has an $\circ$;
    \item $j^\natural$ means the column with $j$ $\times$'s does not have an $\circ$ and is not greyed out; and
    \item $j^\sharp$ means the column with $j$ $\times$'s is greyed out.
  \end{enumerate}
\end{notation}

\begin{reptheorem}{thm:generic_Pieri}[``Generic Pieri Rule'']
  Suppose $I = \{i_1 < \cdots < i_k\} \subset [n-1]$.
  Then
    \[ \sum (-1)^d \partial_{e_{n-k-d}(\underline{n-1})} \dif_{j_1 \cdots j_k} \Delta_n = 0, \]
  where the sum is over all subsets $J = \{j_1 < \cdots < j_k\} \subset [n-1]$
  for which
    \[ 1 \leq i_1 \leq j_1 < i_2 \leq j_2 < i_3 \leq j_3 < \cdots < i_k \leq j_k < n, \]
  where
    \[ d \coloneqq (j_1 - i_1) + \cdots + (j_k - i_k). \]

  \begin{proof}
    We show that for each fixed $0 \leq \ell \leq k$,
    \begin{equation}\label{eq:generic_Pieri.1}
      \sum (-1)^d \{j_1, \ldots, j_\ell, i_{\ell+1}^\natural, \ldots, i_k^\natural\}^{n-k-d} = 0,
    \end{equation}
    where the sum is over $j_1, \ldots, j_\ell$ for which
      \[ 1 \leq i_1 \leq j_1 < \cdots < i_\ell \leq j_\ell < i_{\ell+1}. \]
    Here we define $i_{k+1} \coloneqq n$ and
      \[ d = (j_1 - i_1) + \cdots + (j_\ell - i_\ell). \]
    The theorem is the case $\ell=k$.
    
    We prove \eqref{eq:generic_Pieri.1} by induction on $\ell$. In the base case $\ell=0$, the only term in \eqref{eq:generic_Pieri.1} is $\{i_1^\natural, \ldots, i_k^\natural\}^{n-k}$. Such a marked staircase has $k$ columns with $\times$'s but no $\circ$'s, and $n-k$ columns with $\circ$'s, so at least $n$ columns would have marks. Since an $n$-staircase must have an empty column, there are no such staircases. Now take $\ell \geq 1$. Set $j_0 \coloneqq 0$ for convenience.
    
    We expand each term in \eqref{eq:generic_Pieri.1} using $j_\ell^\flat$, $j_\ell^\natural$, or $j_\ell^\sharp$. Since  $j_{\ell-1} < i_\ell \leq  j_\ell < i_{\ell+1}$, we have $j_\ell = i_\ell, i_\ell+1, \ldots, i_{\ell+1}-1$. If $j_{\ell-1} < j_\ell - 1$, we may apply relation (A) to get
    \begin{equation}\label{eq:generic_Pieri.2}
      \{j_1, \ldots, j_\ell^\natural, i_{\ell+1}^\natural, \ldots, i_k^\natural\}^{n-k-d}
          = \{j_1, \ldots, (j_\ell-1)^\flat, i_{\ell+1}^\natural, \ldots, i_k^\natural\}^{n-k-(d-1)}.
    \end{equation}
    The parity of $(-1)^d$ is opposite for these terms, so they cancel. This observation applies in particular for $j_\ell > i_\ell$. Thus all the terms with $j_\ell^\flat$ or $j_\ell^\natural$ cancel using \eqref{eq:generic_Pieri.2} except for $i_\ell^\natural$ and $(i_{\ell+1}-1)^\flat$. In all, the following terms remain.

    \begin{enumerate}[I.]
      \item $\{j_1, \ldots, j_{\ell-1}, i_\ell^\natural, i_{\ell+1}^\natural, \ldots, i_k^\natural\}^{n-k-d}$. These contributions are $0$ by induction.
      \item $\{j_1, \ldots, j_{\ell-1}, (i_{\ell+1}-1)^\flat, i_{\ell+1}^\natural, \ldots, i_k^\natural\}^{n-k-d}$. When $\ell=k$, we have $i_{k+1}=n$, and $(n-1)^\flat$ would require a column of length $n$, which is too long. For $\ell<k$, we may apply relation (C) to the columns with $i_{\ell+1}-1$ and $i_{\ell+1}$ $\times$'s, which is a sign-reversing involution.
      \item $\{j_1, \ldots, j_{\ell-1}, j_\ell^\sharp, i_{\ell+1}^\natural, \ldots, i_k^\natural\}^{n-k-d}$. We will show that each of these terms is zero. Let $m$ denote the minimum height of the columns with $j_\ell, i_{\ell+1}, \ldots, i_k$ $\times$'s. Since $j_\ell < i_{\ell+1} < \cdots < i_k$, we have $m \geq j_\ell$.
      
        Let $R$ denote the set of columns of height $m, m+1, \ldots, n-1$ which contain an $\circ$. By assumption, the columns with $j_\ell, i_{\ell+1}, \ldots, i_k$ $\times$'s do not contain $\circ$'s, but they would otherwise belong to $R$, so $|R| \leq n-m-(k-\ell+1) = n-k-m+\ell-1$. Since there are $n-k-d$ $\circ$'s,
        \begin{align*}
          \text{$\#\circ$'s} - \# R
            &\geq (n-k - (j_1 - i_1) - \cdots - (j_\ell - i_\ell)) - (n-k-m+\ell-1) \\
            &= (m - j_\ell) + (i_\ell - j_{\ell-1}) + \cdots + (i_2 - j_1) + i_1 - \ell + 1 \\
            &\geq 0 + 1 + \cdots + 1 - \ell + 1 \\
            &= \ell - \ell + 1 \\
            &> 0.
        \end{align*}
        Consequently, there is at least one $\circ$ outside of $R$.
        
        Let $v$ denote the height of the shortest column with an $\circ$. We have just shown $v < m$. By minimality, the column of height $v-1$ (which may be zero) has no $\circ$. Since we have $j_\ell^\sharp$ and there is a unique grey column, the grey column has height at least $j_\ell \geq m > v$, so the columns of height $v$ and $v-1$ are not grey. Thus we may apply relation (B) to swap the $\circ$ between the column of height $v$ and the column of height $v-1$, which is a sign-reversing involution.
    \end{enumerate}
  \end{proof}
\end{reptheorem}

As an application of the Generic Pieri Rule, we prove \eqref{eq:main_cor} for $1$-forms. The relevant case of \eqref{eq:SHn_sgn} is originally due to Alfano \cite{MR1661359}.

\begin{repcorollary}{cor:1_forms}
  The order $\{n-1\} < \{n-2\} < \cdots < \{1\}$ gives a filtration of $\cSH_n^1$ by $\cSH_{\{i\}}$'s where the composition factors are annihilated precisely by the Tanisaki ideal $\cI_{(2, 1^{n-2})}$. In particular,
  \begin{equation*}
    \GrFrob\left(\cSH_n^1; q\right) = [n-1]_q \omega Q_{(2, 1^{n-2})}'(\mathbf{x}; q).
  \end{equation*}
  
  \begin{proof}
    The weakly decreasing rearrangement of $\Phi_n(\{i\})$ is $(2, 1^{n-2})$. The only essential Tanisaki generator for $\cI_{(2, 1^{n-2})}$ is $e_{n-1}(\underline{n-1})$. We have
      \[ \binom{n}{2} - i - b(2, 1^{n-2}) = \binom{n}{2} - i - \binom{n-1}{2} = n-1-i. \]
    Hence the right-hand side of \eqref{eq:main_cor} is
      \[ \sum_{i=1}^{n-1} q^{n-1-i} \omega Q_{(2, 1^{n-2})}'(\mathbf{x}; q) = [n-1]_q \omega Q_{(2, 1^{n-2})}'(\mathbf{x}; q). \]
    
    The minimal term $J$ included in the Generic Pieri Rule is at $J = I = \{i\}$ and is $\partial_{e_{n-1}(\underline{n-1})} \dif_i \Delta_n$. Since all remaining terms have appeared earlier in the filtration, the essential Tanisaki generator annihilates the composition factor. Equality holds in \eqref{eq:SHn_sgn} by Alfano's main result in \cite{MR1661359}, which is equivalent to the $y$-degree $1$ case of the Operator Conjecture/Operator Theorem of Haiman (\cite[Conj.~5.1.1]{MR1256101}, \cite[Thm.~4.2]{MR1918676}) as well as the $\theta$-degree $1$ case of the super operator theorem of Rhoades--Wilson \cite{RW23}. Hence the left-hand side of \eqref{eq:cor:1_forms} is coefficient-wise $\leq$ the right-hand side, as a power series over $q$ in the Schur basis. Alfano in fact showed $\dim \cSH_n^1 = (n-1) n!/2$. Since $\dim \cR_\mu = \binom{n}{\mu}$ and $\binom{n}{2, 1^{n-2}} = n!/2$, equality must hold in \eqref{eq:cor:1_forms}, and the annihilators are tight.
  \end{proof}
\end{repcorollary}

\section{Some symmetric group actions and a shifted Vandermonde identity}\label{sec:actions}

Our proof of the more specific family of extreme hook relations, \Cref{thm:hook_witness_extreme}, is broadly similar to our proof of the Generic Pieri Rule, though it involves grouping certain terms in significantly more intricate ways using certain $\fS_s$-actions and families of involutions. We develop these additional tools now.

\subsection{A shifted Vandermonde identity}

Our upcoming argument will replace a portion of the sets $J \subset [n-1]$ with ordered multisets $\Gamma = (\gamma_1, \ldots, \gamma_s) \subset \mathbb{Z}^s$. We now introduce a family of symmetric group actions on ordered multisets and develop a corresponding shifted Vandermonde evaluation identity, \Cref{cor:shifted_Vandermonde}.

\begin{definition}\label{def:shift_action}
  Suppose $\Gamma = (\gamma_1, \ldots, \gamma_s), \alpha = (\alpha_1, \ldots, \alpha_s) \in \bZ^s$ and $\sigma \in \fS_s$. Define
  \begin{align*}
    \sigma \cdot (\gamma_1, \ldots, \gamma_s) \coloneqq (\gamma_{\sigma^{-1}(1)}, \ldots,
    \gamma_{\sigma^{-1}(s)})
  \end{align*}
  and
  \begin{align}
    \sigma \cdot^\alpha \Gamma \coloneqq \sigma \cdot (\Gamma + \alpha) - \alpha,
  \end{align}
  or explicitly
    \[ \sigma \cdot^\alpha \Gamma = (\gamma_{\sigma^{-1}(1)} + \alpha_{\sigma^{-1}(1)} - \alpha_1,
        \ldots, \gamma_{\sigma^{-1}(s)} + \alpha_{\sigma^{-1}(s)} - \alpha_s). \]
\end{definition}

One may check $\tau \cdot^\alpha (\sigma \cdot^\alpha \Gamma) = (\tau \sigma) \cdot^\alpha \Gamma$, and clearly $\id \cdot^\alpha \Gamma = \Gamma$, so this is a genuine $\fS_s$-action for each fixed $\alpha$. The action $\sigma \cdot^\alpha \Gamma$ is reminiscent of certain actions on weights from Lie theory, e.g. \cite[Cor.~23.2, p.129]{MR0323842}.

\begin{example}
  The $\cdot^\alpha$-orbit of $\Gamma = (2, 2, 3)$ when $\alpha = (1, -1, 0)$ is $\{(2, 2, 3), (2, 4, 1), (0, 4, 3)\}$. Here the stabilizers have order $2$.
\end{example}

\begin{lemma}\label{lem:shifted_Vandermonde}
  Suppose $\Gamma = (\gamma_1, \ldots, \gamma_s), \alpha = (\alpha_1, \ldots, \alpha_s) \in \bZ^s$, and $u \in \bZ_{\geq 0}$. For any fixed $\Pi \subset [s]$ with $|\Pi| > u$,
  \begin{equation}\label{eq:shifted_Vandermonde.1}
    \sum_{\substack{\sigma \in \fS_s \\ M \subset \Pi}} (-1)^{|M|} \sgn(\sigma)
      \Delta_s(\sigma \cdot^\alpha \Gamma - 1_M) |M|^u = 0,
  \end{equation}
  where $\Delta_s(\Gamma) \coloneqq \prod_{1 \leq v < w \leq s} (\gamma_w - \gamma_v)$ is a Vandermonde determinant and $1_M = (\delta_{1 \in M}, \ldots, \delta_{s \in M})$ is the indicator vector for $M$.
  
  \begin{proof}
    Define auxiliary variables $y = (y_1, \ldots, y_s)$ where $y_i \coloneqq \gamma_i + \alpha_i$. Consequently, $\sigma \cdot^\alpha \Gamma = \sigma \cdot y - \alpha$. Now consider the left-hand side of \eqref{eq:shifted_Vandermonde.1} as an element of $\bC[y_1, \ldots, y_s, \alpha_1, \ldots, \alpha_s]$. The $\fS_s$-actions on the $y$ and $\alpha$ variables given by $\tau \circ y_i \coloneqq y_{\sigma(i)}$ and $\rho \circ \alpha_i \coloneqq \alpha_{\rho(i)}$ induce an $\fS_s \times \fS_s$-action on $\bC[y_1, \ldots, y_s, \alpha_1, \ldots, \alpha_s]$.
    
    For $(\tau, \rho) \in \fS_s \times \fS_s$, we have
    \begin{align*}
      (\tau, \rho) \circ \Delta_s(\sigma \cdot^\alpha \Gamma - 1_M)
        &= (\tau, \rho) \circ \Delta_s(\sigma \cdot y - \alpha - 1_M) \\
        &= \Delta_s(\sigma \cdot (\tau^{-1} \cdot y) - \rho^{-1} \cdot \alpha - 1_M) \\
        &= \Delta_s(\rho^{-1} \cdot (\rho \sigma \tau^{-1} \cdot y - \alpha - \rho \cdot 1_M)) \\
        &= \sgn(\rho) \Delta_s(\rho \sigma \tau^{-1} \cdot^\alpha \Gamma - 1_{\rho(M)}),
    \end{align*}
    where in the last line we have used the facts
    \begin{align*}
      \Delta_s(\rho \cdot \Gamma)
        &= \sgn(\rho) \Delta_s(\Gamma), \\
      \rho \cdot 1_M
        &= (\delta_{\rho^{-1}(1) \in M}, \ldots, \delta_{\rho^{-1}(k) \in M}) \\
        &= (\delta_{1 \in \rho(M)}, \ldots, \delta_{k \in \rho(M)}) \\
        &= 1_{\rho(M)}.
    \end{align*}
    Consequently,
    \begin{align*}
      (\tau, \rho) &\circ \sum_{\substack{\sigma \in \fS_s \\ M \subset \Pi}}
          (-1)^{|M|} \sgn(\sigma) \Delta_s(\sigma \cdot^\alpha \Gamma - 1_M) |M|^u \\
        &= \sum_{\substack{\sigma \in \fS_s \\ M \subset \Pi}}
          (-1)^{|M|} \sgn(\rho)\sgn(\sigma)
          \Delta_s(\rho \sigma \tau^{-1} \cdot^\alpha \Gamma - 1_{\rho(M)}) |M|^u \\
        &= \sum_{\substack{\sigma \in \fS_s \\ M \subset \rho(\Pi)}}
          (-1)^{|M|} \sgn(\rho)\sgn(\rho^{-1}\sigma\tau)
          \Delta_s(\sigma \cdot^\alpha \Gamma - 1_M) |M|^u \\
        &= \sgn(\tau) \sum_{\substack{\sigma \in \fS_s \\ M \subset \rho(\Pi)}}
          (-1)^{|M|} \sgn(\sigma)
          \Delta_s(\sigma \cdot^\alpha \Gamma - 1_M) |M|^u
    \end{align*}
    where in the second step we have reindexed according to
    $M \mapsto \rho^{-1}(M)$ and $\sigma \mapsto \rho^{-1}\sigma\tau$.
    
    Letting $\rho=\id$ and specializing the $\alpha$ variables to integer constants, this last expression says the left-hand side of \eqref{eq:shifted_Vandermonde.1} as an inhomogeneous element of $\bC[y_1, \ldots, y_s]$ is an alternating polynomial. Thus all components of $y$-degree below $\deg \Delta_s(y_1, \ldots, y_s)$ vanish. The only possible remaining component is
    \begin{align*}
      \sum_{\substack{\sigma \in \fS_s \\ M \subset \Pi}} (-1)^{|M|} \sgn(\sigma)
        \Delta_s(\sigma \cdot y) |M|^u
        &= \sum_{\substack{\sigma \in \fS_s \\ M \subset \Pi}} (-1)^{|M|} \Delta_s(y) |M|^u \\
        &= s! \Delta_s(y) \sum_{M \subset \Pi} (-1)^{|M|} |M|^u.
    \end{align*}
    It is well-known that $\sum_{M \subset \Pi} (-1)^{|M|} |M|^u = 0$ for $|\Pi| > u$. Indeed, it is $(-1)^{p} p! \Stir(u, p)$ where $\Stir$ denotes a Stirling number of the second kind and $p \coloneqq |\Pi|$. More directly, it follows from differentiating the binomial theorem
      \[ (1+x)^p = \sum_{k=0}^{p} \binom{p}{k} x^k \]
    up to $u$ times and setting $x=-1$.
  \end{proof}
\end{lemma}

\begin{corollary}\label{cor:shifted_Vandermonde}
  For any $\Pi \subset [s]$ with $|\Pi| > u \geq 0$, $v \in \mathbb{Z}$, and $\Gamma, \alpha \in \mathbb{Z}^s$,
    \[ \sum_{\substack{\sigma \in \fS_s \\ M \subset \Pi}} (-1)^{|M|} \sgn(\sigma)
        \Delta_s(\sigma \cdot^\alpha \Gamma - 1_M) \binom{v - |M| + u}{u}
        = 0. \]
  
  \begin{proof}
    The factor $\binom{v - |M| + u}{u} = \frac{1}{u!} \prod_{i=1}^u (v - |M| + i)$ is a polynomial in $|M|$ of degree $u < |\Pi|$. The result follows by taking linear combinations of \Cref{lem:shifted_Vandermonde}.
  \end{proof}
\end{corollary}

\subsection{An action on marked staircases}

Our upcoming argument will group together certain marked staircases using another $\fS_s$-action. We introduce this action with the following technical lemma. See \Cref{ex:staircase_action} and \Cref{fig:staircase_action}.

\begin{lemma}\label{lem:staircase_action}
  Suppose $I = \{i_1 < \cdots < i_k\} \subset [n-1]$ where $\alpha = \Phi_n(I)$ satisfies $\overline{\alpha} = (s, 1^{k-s})$ for some $1 \leq s \leq k$.
  
  Let $M_I$ denote the set of all marked staircases where the multiset $J = \{\!\{j_1 \leq \cdots \leq j_k\}\!\}$ of the number of $\times$'s in each column satisfies
    \[ j_1 = i_1, \ldots, j_{k-s} = i_{k-s}, \]
  where
    \[ d \coloneqq \sum_{\ell=k-s+1}^k (j_\ell - i_\ell) \geq 0. \]
  Then:
  \begin{enumerate}[(i)]
    \item Every realizable multiset $\{\!\{j_1 \leq \cdots \leq j_k\}\!\}$ is lexicographically greater than or equal to the set $\{i_{k-s+1} < \cdots < i_k\}$. More precisely, if $d>0$, then $j_{k-s+1} > i_{k-s+1}$, and if $d=0$, then $I = J$.
    \item The unique set of columns with $j_{k-s+1}, \ldots, j_k$ $\times$'s all have some $\times$ at the same, common height.
    \item Relation (D) gives an $\fS_s$-action on $M_I$ by acting on the columns with $j_{k-s+1}, \ldots, j_k$ $\times$'s.
    \item Moreover, $d < m$ where $m$ is the minimum height of a column with $j_{k-s+1}, \ldots, j_k$ $\times$'s.
  \end{enumerate}
  
  \begin{proof}
    We may \textit{a priori} have $j_{k-s+1} = j_{k-s}$, in which case the set of columns from (ii) and (iii) is not unique. For now, choose some set of columns with $j_{k-s+1}, \ldots, j_k$ $\times$'s and let the set of heights of these columns be $\{h_1 < \cdots < h_s\}$.
    
    First consider $s=1$. Here $d \geq 0$ gives $j_{k-s+1} \geq i_{k-s+1}$, so (i) holds. Since $i_{k-s+1} > i_{k-s}$, uniqueness holds in (ii) and the remaining conclusions in (ii), (iii), and (iv) are trivial or obvious. Now suppose $s \geq 2$.
    
    Since the marks fit in an $n$-staircase, we have $h_1 \leq n-s, h_2 \leq n-s+1, \ldots, h_s \leq n-1$. Write $j_\ell'$ for the number of $\times$'s in the column of height $h_\ell$, so $j_\ell' \leq h_\ell \leq n-s+\ell-1$. Hence we have $\delta_2, \ldots, \delta_s \geq 0$ and some $\epsilon_1 \in \mathbb{Z}$ for which
      \[ (j_1', j_2', \ldots, j_s') = (i_{k-s+1} + \epsilon_1,
           n-s+1-\delta_2, \ldots, n-1-\delta_s). \]
    
    Recall from \Cref{rem:extreme_hook} that the condition on $I$ forces $i_{k-s+2} = n-s+1, \ldots, i_k = n-1$, so $i_{k-s+\ell} = n-s+\ell-1$ for $2 \leq \ell \leq s$. Hence we have
    \begin{align*}
      \delta_\ell &= i_{k-s+\ell} - j_\ell' \qquad (2 \leq \ell \leq s) \\
      \epsilon_1 &= j_1' - i_{k-s+1}.
    \end{align*}
    Thus
      \[ d = \epsilon_1 - \delta_2 - \cdots - \delta_s \geq 0, \]
    so $\epsilon_1 \geq \delta_2 + \cdots + \delta_s \geq 0$. In particular, $\epsilon_1 \geq \delta_\ell$ for all $2 \leq \ell \leq s$ and $j_1' \geq i_{k-s+1}$. Consequently,
    \begin{align*}
      d &\leq \epsilon_1 = j_1' - i_{k-s+1} < j_1' \leq h_1 = m,
    \end{align*}
    giving (iv), assuming uniqueness for the moment.
    
    For a fixed $2 \leq \ell \leq s$, consider the height of the lowest $\times$ in the column with height $h_\ell$ and $j_\ell'$ $\times$'s. This height is
    \begin{align}
      h_\ell - j_\ell' + 1 \nonumber
        &\leq i_{k-s+\ell} - j_\ell' + 1 \nonumber \\
        &= \delta_\ell + 1 \leq \epsilon_1 + 1 \nonumber \\
        &= j_1' - i_{k-s+1} + 1 \nonumber \\
        &\leq h_1 - i_{k-s+1} + 1. \label{eq:staircase_action.1}
    \end{align}
    On the other hand, the highest $\times$ in this column is at the top at height $h_\ell \geq h_1$. Hence we must at least have $\times$'s in this column at heights $h_1 - i_{k-s+1} + 1, \ldots, h_1$, or at least $i_{k-s+1}$ of them in all. The same is true of the column with height $h_1$ since $j_1' \geq i_{k-s+1}$, so all columns of height $h_1, \ldots, h_s$ have $\times$'s at these $i_{k-s+1}$ common heights. This proves proves (ii), except for the uniqueness claim.
    
    The preceding argument gives $j_\ell' \geq i_{k-s+1}$ for $1 \leq \ell \leq s$. Indeed, since $h_1 < \cdots < h_\ell$, \eqref{eq:staircase_action.1} gives the tighter bound
    \begin{align*}
      j_\ell'
        &\geq h_\ell - h_1 + i_{k-s+1} \\
        &\geq i_{k-s+1} + (\ell-1),
    \end{align*}
    which also holds at $\ell=1$. Hence
    \begin{align*}
      j_{k-s+1}
        &= \min\{j_{k-s+\ell} : 1 \leq \ell \leq s\} \\
        &= \min\{j_\ell' : 1 \leq \ell \leq s \} \\
        &\geq i_{k-s+1} > i_{k-s} = j_{k-s},
    \end{align*} 
    so uniqueness follows as well and (ii) holds.
    
    Moreover, since $j_\ell' \geq i_{k-s+1} + (\ell-1)$, we see that $\{\!\{j_1', \ldots, j_s'\}\!\} = \{\!\{j_{k-s+1} \leq \cdots \leq j_k\}\!\}$ is lexicographically larger than $\{i_{k-s+1} < \cdots < i_k\}$ except perhaps when $j_1' = j_{k-s+1} = i_{k-s+1}$. However, in that case $\epsilon_1 = 0$, forcing $\delta_\ell = 0$ for $2 \leq \ell \leq s$, so $j_\ell' = i_{k-s+\ell}$. Thus (i) holds.

    Finally, by (ii) we may take the $s$ columns of heights $h_1, \ldots, h_s$ and permute them amongst themselves using relation (D), giving an $\fS_s$-action; see \Cref{ex:staircase_action} and \Cref{fig:staircase_action}. We must only show that the resulting marked staircase remains in $M_I$. The action preserves $d$ and $j_1, \ldots, j_{n-k}$, so we must only show that among the permuted columns, none have fewer than $j_{n-k}$ $\times$'s. But we showed above that the columns have a block of $i_{n-k+1} > i_{n-k} = j_{n-k}$ $\times$'s at a common height, which is preserved by relation (D), giving (iii) and completing the proof.
  \end{proof}
\end{lemma}

\begin{definition}\label{def:staircase_action}
  Suppose $S \in M_I$ from \Cref{lem:staircase_action}. Define an explicit $\fS_s$-action on $M_I$ as follows.
  \begin{itemize}
    \item Call the $s$ columns with at least $i_{k-s+1}$ $\times$'s the \textit{active columns}.
    \item Let $h_\ell$ be the height of the $\ell$th active column. Let $\sigma \cdot S$ be the marked staircase obtained by applying relation (D) to $S$ where the $\sigma(\ell)$th active column has height $h_\ell$.
  \end{itemize}
  Furthermore:
  \begin{itemize}
    \item Let $\Gamma(S) = (\gamma_1, \ldots, \gamma_s)$ be the number of $\times$'s in the active columns of $S$, read from left to right. Note that $(j_{k-s+1}, \ldots, j_k)$ is the weakly increasing rearrangement of $\Gamma$.
    \item Let $\alpha(S) = (\alpha_1, \ldots, \alpha_s)$ be the number of cells without $\times$'s in the active columns of $S$, read from left to right. Note that $\alpha(\sigma \cdot S) = \alpha(S)$.
  \end{itemize}
\end{definition}

\begin{example}\label{ex:staircase_action}
  Let $I = \{1, 3, 4, 8, 9\} \subset [11-1]$, so $n=10$ and $k=5$. We have $\Phi_n(I) = (1, 4, 2, 1, 2)$, $\overline{\Phi_n(I)} = (3, 1, 1)$, and $s=3$. The corresponding diagram is:
  \begin{align*}
    \tableau{
      \ \\
      \ & 4 & 8 & 9 \\
      \ & 3 \\
      \ \\
      \ & 1}
  \end{align*}
  See \Cref{fig:staircase_action} for an $\fS_3$-orbit of marked staircases in $M_I$ from \Cref{lem:staircase_action}.
\begin{figure}[ht]
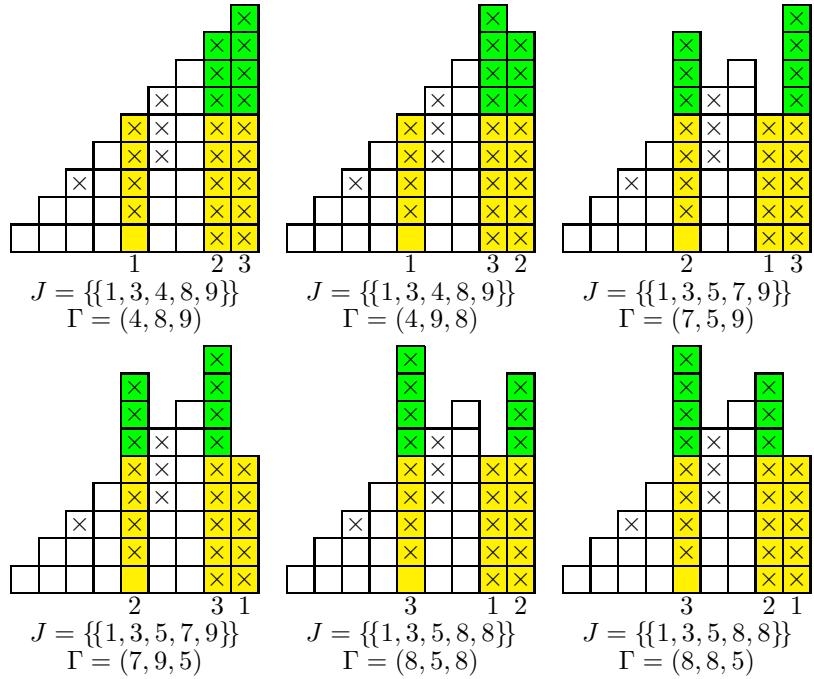

  \centering

  \begin{align*}
    \begin{ytableau}
      \none & \none & \none & \none & \none & \none & \none & \none & \none & *(green) \times \\
      \none & \none & \none & \none & \none & \none & \none & \none & *(green) \times & *(green) \times \\
      \none & \none & \none & \none & \none & \none & \none & \  & *(green) \times & *(green) \times \\
      \none & \none & \none & \none & \none & \none & \times & \  & *(green) \times & *(green) \times \\
      \none & \none & \none & \none & \none & *(yellow) \times & \times & \  & *(yellow) \times & *(yellow) \times \\
      \none & \none & \none & \none & \  & *(yellow) \times &  \times & \  & *(yellow) \times & *(yellow) \times \\
      \none & \none & \none & \times & \  & *(yellow) \times & \  & \  & *(yellow) \times & *(yellow) \times \\
      \none & \none & \  & \  & \  & *(yellow) \times & \  & \  & *(yellow) \times & *(yellow) \times \\
      \none & \  & \  & \  & \  & *(yellow) \  & \  & \  & *(yellow) \times & *(yellow) \times \\
      \none & \none & \none & \none & \none & \none[1] & \none & \none & \none[2] & \none[3] \\
      \none & \none & \none & \none & \none & \none[J = \{\!\{1, 3, 4, 8, 9\}\!\}] & \none & \none & \none & \none \\
      \none & \none & \none & \none & \none & \none[\Gamma = (4, 8, 9)] & \none & \none & \none & \none
    \end{ytableau}
    \begin{ytableau}
      \none & \none & \none & \none & \none & \none & \none & \none & *(green) \times & \none \\
      \none & \none & \none & \none & \none & \none & \none & \none & *(green) \times & *(green) \times \\
      \none & \none & \none & \none & \none & \none & \none & \  & *(green) \times & *(green) \times \\
      \none & \none & \none & \none & \none & \none & \times & \  & *(green) \times & *(green) \times \\
      \none & \none & \none & \none & \none & *(yellow) \times & \times & \  & *(yellow) \times & *(yellow) \times \\
      \none & \none & \none & \none & \  & *(yellow) \times &  \times & \  & *(yellow) \times & *(yellow) \times \\
      \none & \none & \none & \times & \  & *(yellow) \times & \  & \  & *(yellow) \times & *(yellow) \times \\
      \none & \none & \  & \  & \  & *(yellow) \times & \  & \  & *(yellow) \times & *(yellow) \times \\
      \none & \  & \  & \  & \  & *(yellow) \  & \  & \  & *(yellow) \times & *(yellow) \times \\
      \none & \none & \none & \none & \none & \none[1] & \none & \none & \none[3] & \none[2] \\
      \none & \none & \none & \none & \none & \none[J = \{\!\{1, 3, 4, 8, 9\}\!\}] & \none & \none & \none & \none \\
      \none & \none & \none & \none & \none & \none[\Gamma = (4, 9, 8)] & \none & \none & \none & \none
    \end{ytableau}
    \begin{ytableau}
      \none & \none & \none & \none & \none & \none & \none & \none & \none & *(green) \times \\
      \none & \none & \none & \none & \none & *(green) \times & \none & \none & \none & *(green) \times \\
      \none & \none & \none & \none & \none & *(green) \times & \none & \  & \none & *(green) \times \\
      \none & \none & \none & \none & \none & *(green) \times & \times & \  & \none & *(green) \times \\
      \none & \none & \none & \none & \none & *(yellow) \times & \times & \  & *(yellow) \times & *(yellow) \times \\
      \none & \none & \none & \none & \  & *(yellow) \times &  \times & \  & *(yellow) \times & *(yellow) \times \\
      \none & \none & \none & \times & \  & *(yellow) \times & \  & \  & *(yellow) \times & *(yellow) \times \\
      \none & \none & \  & \  & \  & *(yellow) \times & \  & \  & *(yellow) \times & *(yellow) \times \\
      \none & \  & \  & \  & \  & *(yellow) \  & \  & \  & *(yellow) \times & *(yellow) \times \\
      \none & \none & \none & \none & \none & \none[2] & \none & \none & \none[1] & \none[3] \\
      \none & \none & \none & \none & \none & \none[J = \{\!\{1, 3, 5, 7, 9\}\!\}] & \none & \none & \none & \none \\
      \none & \none & \none & \none & \none & \none[\Gamma = (7, 5, 9)] & \none & \none & \none & \none
    \end{ytableau} \\
    \begin{ytableau}
      \none & \none & \none & \none & \none & \none & \none & \none & *(green) \times & \none \\
      \none & \none & \none & \none & \none & *(green) \times & \none & \none & *(green) \times & \none \\
      \none & \none & \none & \none & \none & *(green) \times & \none & \  & *(green) \times & \none \\
      \none & \none & \none & \none & \none & *(green) \times & \times & \  & *(green) \times & \none \\
      \none & \none & \none & \none & \none & *(yellow) \times & \times & \  & *(yellow) \times & *(yellow) \times \\
      \none & \none & \none & \none & \  & *(yellow) \times &  \times & \  & *(yellow) \times & *(yellow) \times \\
      \none & \none & \none & \times & \  & *(yellow) \times & \  & \  & *(yellow) \times & *(yellow) \times \\
      \none & \none & \  & \  & \  & *(yellow) \times & \  & \  & *(yellow) \times & *(yellow) \times \\
      \none & \  & \  & \  & \  & *(yellow) \  & \  & \  & *(yellow) \times & *(yellow) \times \\
      \none & \none & \none & \none & \none & \none[2] & \none & \none & \none[3] & \none[1] \\
      \none & \none & \none & \none & \none & \none[J = \{\!\{1, 3, 5, 7, 9\}\!\}] & \none & \none & \none & \none \\
      \none & \none & \none & \none & \none & \none[\Gamma = (7, 9, 5)] & \none & \none & \none & \none
    \end{ytableau}
    \begin{ytableau}
      \none & \none & \none & \none & \none & *(green) \times & \none & \none & \none & \none \\
      \none & \none & \none & \none & \none & *(green) \times & \none & \none & \none & *(green) \times \\
      \none & \none & \none & \none & \none & *(green) \times & \none & \  & \none & *(green) \times \\
      \none & \none & \none & \none & \none & *(green) \times & \times & \  & \none & *(green) \times \\
      \none & \none & \none & \none & \none & *(yellow) \times & \times & \  & *(yellow) \times & *(yellow) \times \\
      \none & \none & \none & \none & \  & *(yellow) \times &  \times & \  & *(yellow) \times & *(yellow) \times \\
      \none & \none & \none & \times & \  & *(yellow) \times & \  & \  & *(yellow) \times & *(yellow) \times \\
      \none & \none & \  & \  & \  & *(yellow) \times & \  & \  & *(yellow) \times & *(yellow) \times \\
      \none & \  & \  & \  & \  & *(yellow) \  & \  & \  & *(yellow) \times & *(yellow) \times \\
      \none & \none & \none & \none & \none & \none[3] & \none & \none & \none[1] & \none[2] \\
      \none & \none & \none & \none & \none & \none[J = \{\!\{1, 3, 5, 8, 8\}\!\}] & \none & \none & \none & \none \\
      \none & \none & \none & \none & \none & \none[\Gamma = (8, 5, 8)] & \none & \none & \none & \none
    \end{ytableau}
    \begin{ytableau}
      \none & \none & \none & \none & \none & *(green) \times & \none & \none & \none & \none \\
      \none & \none & \none & \none & \none & *(green) \times & \none & \none & *(green) \times & \none \\
      \none & \none & \none & \none & \none & *(green) \times & \none & \  & *(green) \times & \none \\
      \none & \none & \none & \none & \none & *(green) \times & \times & \  & *(green) \times & \none \\
      \none & \none & \none & \none & \none & *(yellow) \times & \times & \  & *(yellow) \times & *(yellow) \times \\
      \none & \none & \none & \none & \  & *(yellow) \times &  \times & \  & *(yellow) \times & *(yellow) \times \\
      \none & \none & \none & \times & \  & *(yellow) \times & \  & \  & *(yellow) \times & *(yellow) \times \\
      \none & \none & \  & \  & \  & *(yellow) \times & \  & \  & *(yellow) \times & *(yellow) \times \\
      \none & \  & \  & \  & \  & *(yellow) \  & \  & \  & *(yellow) \times & *(yellow) \times \\
      \none & \none & \none & \none & \none & \none[3] & \none & \none & \none[2] & \none[1] \\
      \none & \none & \none & \none & \none & \none[J = \{\!\{1, 3, 5, 8, 8\}\!\}] & \none & \none & \none & \none \\
      \none & \none & \none & \none & \none & \none[\Gamma = (8, 8, 5)] & \none & \none & \none & \none
    \end{ytableau}
  \end{align*}
  \caption{An $\fS_3$-orbit in $M_I$ for $I$ from \Cref{ex:staircase_action} obtained by applying relation (D) as in \Cref{lem:staircase_action} to the three active columns with the most $\times$'s in each marked staircase. The relative order of the three active columns forms a permutation $\sigma^{-1}$ which has been written below the marked staircases. The staircases are of the form $\sigma \cdot S$ where $S$ is the upper-left diagram. The multiset of the number of $\times$'s in all columns is $J$ and the number of $\times$'s in the three active columns from left to right is $\Gamma$. In each case, the number of cells without $\times$'s in the active columns is $\alpha = (1, 0, 0)$.}
  \label{fig:staircase_action}
\end{figure}
\end{example}

The actions from \Cref{def:shift_action} and \Cref{lem:staircase_action} are related as follow. See \Cref{ex:staircase_action.2}.

\begin{lemma}\label{lem:staircase_action.2}
  Let $S \in M_I$ from \Cref{lem:staircase_action} and $\sigma \in \fS_n$ with $\alpha = \alpha(S) = \alpha(\sigma \cdot S)$. Then
  \begin{align*}
    \Gamma(\sigma \cdot S)
      &= \sigma \cdot^\alpha \Gamma(S)
    \end{align*}
    and
    \begin{align*}
      \sgn \Delta_s(\Gamma(S)) \wgt(\sigma \cdot S) &= \sgn(\sigma) \sgn \Delta_s(\sigma \cdot^\alpha \Gamma(S)) \wgt(S).
  \end{align*}

  \begin{proof}
    Let $\Gamma(S) = (\gamma_1, \ldots, \gamma_s)$ and $\Gamma(\sigma \cdot S) = (\lambda_1, \ldots, \lambda_s)$. We have $\alpha_\ell + \gamma_\ell = h_\ell$ where $h_\ell$ is the height of the $\ell$th active column of $S$. Similarly, we have $\alpha_{\sigma(\ell)} + \lambda_{\sigma(\ell)} = h_\ell$. Hence
    \begin{align*}
      \lambda_\ell
        &= h_{\sigma^{-1}(\ell)} - \alpha_\ell \\
        &= \gamma_{\sigma^{-1}(\ell)} + \alpha_{\sigma^{-1}(\ell)} - \alpha_\ell \\
        &= (\sigma \cdot (\Gamma(S) + \alpha))_\ell - \alpha_\ell \\
        &= (\sigma \cdot^\alpha \Gamma(S))_\ell,
    \end{align*}
    giving the first claim.
    
    For the second statement, let $N$ be the order of $S$ and $\sigma \cdot S$ and let $c_1 < \cdots < c_k$ be the indexes of the columns with $\times$'s. Let $j_\ell$ be the number of $\times$'s in column $c_\ell$ of $S$ and let $j_\ell'$ be the number of $\times$'s in column $c_\ell$ of $\sigma \cdot S$. Then we have $c$ and $\beta$ where
    \begin{align*}
      \wgt(S)
        &= (-1)^c Nx^\beta \sgn \Delta_k(j_1, \ldots, j_k) \theta_{c_1} \cdots \theta_{c_k} \\
      \wgt(\sigma \cdot S)
        &= (-1)^c \sgn(\sigma) Nx^\beta \sgn \Delta_k(j_1', \ldots, j_k') \theta_{c_1} \cdots \theta_{c_k}.
    \end{align*}
    If $\ell$ is not an active column, then $j_\ell' = j_\ell < i_{k-s+1}$. If $\ell$ is an active column, then $j_\ell, j_\ell' \geq i_{k-s+1}$. It follows that the sign difference between $\Delta_k(j_1, \ldots, j_k)$ and $\Delta_k(j_1', \ldots, j_k')$ is precisely the same as the sign difference between $\Delta_s(\gamma_1, \ldots, \gamma_s)$ and $\Delta_s(\lambda_1, \ldots, \lambda_s)$. The result follows by combining these observations.
  \end{proof}
\end{lemma}

\begin{example}\label{ex:staircase_action.2}
  Let $S$ be the upper left diagram in \Cref{fig:staircase_action} and $\sigma = 312 = 231^{-1}$, so $\sigma \cdot S$ is the lower left diagram. We have
    \[ (7, 9, 5) = 231 \cdot^{(1, 0, 0)} (4, 8, 9), \]
  in agreement with \Cref{lem:staircase_action.2}.
\end{example}

\section{Extreme hook relations proof}\label{sec:hook}

We may finally prove our second family of Tanisaki witness relations, \Cref{thm:hook_witness_extreme}. The argument will rely on grouping marked staircases using the following more technical variation on \Cref{not:generic_Pieri}.

\begin{notation}\label{not:hook_witness_extreme}
  Fix $I = \{i_1 < \cdots < i_k\} \subset [n-1]$ with $\overline{\Phi_n(I)} = (s, 1^{k-s})$ for some $1 \leq s \leq k$ as in \Cref{lem:staircase_action} and \Cref{def:staircase_action}. Each marked staircase $S \in M_I$ has the following data attached to it.
  \begin{itemize}
    \item The multiset $J = \{\!\{j_1 \leq \cdots \leq j_k\}\!\} \subset [n-1]$ giving the number of $\times$'s in columns with them.
    \item The number $\delta$ of $\circ$'s.
    \item The number $\eta$ of grey columns.
    \item The set of $s$ active columns, namely those with at least $i_{k-s+1}$ $\times$'s.
    \item The list $\Gamma = (\gamma_1, \ldots, \gamma_s)$ of the number of $\times$'s in the $s$ active columns, read from left to right.
    \item The list $\alpha = (\alpha_1, \ldots, \alpha_s)$ of the number of cells in the active columns without $\times$'s, read from left to right.
    \item The subset $\Omega$ of active columns with $\circ$'s.
    \item The subset $\Pi$ of non-grey active columns without $\circ$'s.
    \item The subset $\Psi$ of grey active columns.
  \end{itemize}
  We consider $\Omega \sqcup \Pi \sqcup \Psi = [s]$ by numbering the active columns left to right from $1$ to $s$. Let
    \[ (\Gamma, \Omega, \Pi)_\alpha^{\delta, \eta} \]
  denote the weight generating function of marked staircases in $M_I$ with the above data.
\end{notation}

\begin{reptheorem}{thm:hook_witness_extreme}
  Suppose $I = \{i_1 < \cdots < i_k\} \subset [n-1]$ is such that for some $1 \leq s \leq k$ we have
  \begin{align*}
    i_1, \ldots, i_{k-s+1} &\leq n-k \\
    i_{k-s+2} &= n-s+1 \\
    i_{k-s+3} &= n-s+2 \\
    &\ \,\vdots \\
    i_k &= n-1.
  \end{align*}
  Pick $0 \leq u \leq s$. Then
  \begin{equation*}
    \sum (-1)^d \Delta_s(j_{k-s+1}, \ldots, j_k) \binom{d+u}{u}
        \partial_{e_{n-s-d}(\underline{n-s+u})} \dif_J \Delta_n = 0,
  \end{equation*}
  where the sum is over all subsets $J = \{j_1 < \cdots < j_k\} \subset [n-1]$
  for which
  \begin{align*}
    j_1 = i_1, \ldots, j_{k-s} &= i_{k-s} \\
    d \coloneqq (j_{k-s+1} - i_{k-s+1}) + &\cdots + (j_k - i_k) \geq 0.
  \end{align*}
  
  \begin{proof}
    Since $\dif_J = 0$ if terms repeat, we may include multisets $J = \{\!\{j_1 \leq \cdots \leq j_k\}\!\} \subset [n-1]$ in \eqref{eq:hook_witness_extreme.1}. By \Cref{lem:marked_staircases_gf} and \Cref{lem:staircase_action},
    \begin{equation}\label{eq:hook_witness_extreme.1.1}
      \partial_{e_\delta(\underline{n-\eta})} \dif_J \Delta_n = \sum_{\substack{\Omega, \Pi, \alpha \\ \Gamma \sim J_{\mathrm{top}}}} (\Gamma, \Omega, \Pi)_\alpha^{\delta, \eta},
    \end{equation}
where $J_{\mathrm{top}} \coloneqq (j_{n-k+1}, \ldots, j_k)$ and $\Gamma \sim J_{\mathrm{top}}$ means the weakly increasing rearrangement of $\Gamma$ is $J_{\mathrm{top}}$.

    If $\Gamma \sim J_{\mathrm{top}}$, then
    \begin{equation}\label{eq:hook_witness_extreme.1.2}
      \Delta_s(j_{k-s+1}, \ldots, j_k)
        = \Delta_s(J_{\mathrm{top}}) 
        = \sgn \Delta_s(\Gamma) \cdot \Delta_s(\Gamma).
    \end{equation}
    Using \eqref{eq:hook_witness_extreme.1.1} and \eqref{eq:hook_witness_extreme.1.2}, the left-hand side of \eqref{eq:hook_witness_extreme.1} becomes
    \begin{align}
      \sum_J (-1)^d &\Delta_s(j_{k-s+1}, \ldots, j_k) \binom{d+u}{u} \partial_{e_{n-s-d}(\underline{n-s+u})} \dif_J \Delta_n \nonumber \\
        &= \sum_J (-1)^d \Delta_s(J_\mathrm{top}) \binom{d+u}{u} \sum_{\substack{\Omega, \Pi, \alpha \\ \Gamma \sim J_{\mathrm{top}}}} (\Gamma, \Omega, \Pi)_\alpha^{n-s-d, s-u} \nonumber \\
        &= \sum_{\Gamma, \Omega, \Pi, \alpha} (-1)^d \sgn \Delta_s(\Gamma) \cdot \Delta_s(\Gamma) \binom{d+u}{u} (\Gamma, \Omega, \Pi)_\alpha^{n-s-d, s-u}, \label{eq:hook_witness_extreme.2}
    \end{align}
    where
      \[ d = \Sum(J) - \Sum(I) = \Sum(\Gamma) - \Sum(I_{\mathrm{top}}). \]
    We will group the contributions to \eqref{eq:hook_witness_extreme.2} into terms which individually sum to zero.
    
    As a warm-up, we first show that we may apply relation (B) to cancel all contributions when $\Omega = \varnothing$ and $|\Pi| = u$. In this case, no active columns have $\circ$'s and every grey column is active. Let $m$ denote the minimal height of the active columns. Let $R$ be the set of columns of height $m, m+1, \ldots, n-1$ which contain an $\circ$. All $s$ active columns are in this height range but have no $\circ$'s, so $\#R \leq n-m-s$. There are $n-s-d$ $\circ$'s, so
    \begin{align*}
      \#\circ\text{'s} - \#R
        &\geq (n-s-d) - (n-m-s) \\
        &= m-d \\
        &> 0,
    \end{align*}
    where we have used \Cref{lem:staircase_action}(iv). Thus there are columns with $\circ$'s outside of $R$. Let $v$ be the minimum height of a column with an $\circ$. We have just shown $v < m$, so the columns with heights $v$ and $v-1$ (the latter may have height $0$) are not active, and hence are not grey. We may now apply relation (B) to cancel these terms. That is,
      \[ (\Gamma, \varnothing, \Pi)_\alpha^{n-s-d, s-u} = 0 \qquad\text{if }|\Pi| = u. \]

    The $\fS_s$-action from \Cref{lem:staircase_action} preserves $\Omega, \Pi, d, \alpha$. By \Cref{lem:staircase_action.2}, the action replaces $\Gamma$ with $\sigma \cdot^\alpha \Gamma$. If $\Gamma$ contains repeated elements, then $(\Gamma, \Omega, \Pi)_\alpha^{\delta, \eta} = 0$, so we assume $\Gamma$ does not contain repeated elements. Now \Cref{lem:staircase_action.2} gives
      \[ \wgt(\sigma \cdot S) = \sgn(\sigma) \sgn \Delta_s(\Gamma) \sgn \Delta_s(\sigma \cdot^\alpha \Gamma) \wgt(S). \]
    Hence
    \begin{equation}\label{eq:hook_witness_extreme.3.1}
      (\sigma \cdot^\alpha \Gamma, \Omega, \Pi)_\alpha^{n-s-d, s-u} = \sgn(\sigma) \sgn \Delta_s(\Gamma) \sgn \Delta_s(\sigma \cdot^\alpha \Gamma) (\Gamma, \Omega, \Pi)_\alpha^{n-s-d, s-u}.
    \end{equation}
    
    Now pick some subset $M \subset \Pi$ of the non-grey active columns without $\circ$'s and apply relation (A) to each of those active columns, replacing their bottom-most $\times$'s with an $\circ$.  The resulting staircase remains in $M_I$ so long as $d$ remains non-negative. In this case, the operation replaces $\alpha$ with $\alpha + 1_M$, $\Gamma$ with $\Gamma - 1_M$, $\Pi$ with $\Pi - M$, $\Omega$ with $\Omega \sqcup M$, and $d$ with $d-|M|$. Hence if we require $|M| \leq \Sum(J) - \Sum(I_{\mathrm{top}})$, this operation is well-defined and indeed invertible. This operation preserves monomial weight, and we have
    \begin{align}
      (\Gamma - 1_M, \Omega \sqcup M, &\Pi - M)_{\alpha+1_M}^{n-s-d+|M|, s-u} \nonumber \\
        &= \sgn \Delta_s(\Gamma) \sgn \Delta_s(\Gamma - 1_M) (\Gamma, \Omega, \Pi)_\alpha^{n-s-d, s-u}.\label{eq:hook_witness_extreme.3.2}
    \end{align}
    Combining \eqref{eq:hook_witness_extreme.3.1} and \eqref{eq:hook_witness_extreme.3.2}, we have
    \begin{align}
      (\sigma \cdot^\alpha \Gamma - 1_M, \Omega \sqcup M, &\Pi - M)_{\alpha+1_M}^{n-s-d+|M|, s-u} \nonumber \\
        &= \sgn(\sigma) \sgn \Delta_s(\Gamma) \sgn \Delta_s(\sigma \cdot^\alpha \Gamma - 1_M) (\Gamma, \Omega, \Pi)_\alpha^{n-s-d, s-u}. \label{eq:hook_witness_extreme.3.3}
     \end{align}
    
    Suppose now that $\Omega = \varnothing$ and $|\Pi| > u$. Consider the contributions to \eqref{eq:hook_witness_extreme.2} arising from the ``orbit'' obtained by first applying the $\fS_s$-action and then applying relation (A) as above. By \eqref{eq:hook_witness_extreme.3.3}, these contributions are
    \begin{alignat*}{2}
      &\sum_{\substack{\sigma \in \fS_n \\ M \subset \Pi \\ d-|M| \geq 0}}
        &&(-1)^{d-|M|} \sgn \Delta_s(\sigma \cdot^\alpha \Gamma - 1_M) \cdot \Delta_s(\sigma \cdot^\alpha \Gamma - 1_M) \binom{d-|M|+u}{u} \\
        &&&\cdot (\sigma \cdot^\alpha \Gamma, M, \Pi - M)_{\alpha+1_M}^{n-s-d+|M|, s-u} \\
        =&(-1)^d &&\sgn \Delta_s(\Gamma) \sum_{\substack{\sigma \in \fS_n \\ M \subset \Pi \\ d-|M| \geq 0}}
        (-1)^{|M|} \sgn(\sigma) \Delta_s(\sigma \cdot^\alpha \Gamma - 1_M) \binom{d-|M|+u}{u} \\
        &&& \cdot (\Gamma, \varnothing, \Pi)_\alpha^{n-s-d, s-u}.
    \end{alignat*}
    Here we must interpret the binomial coefficient as the polynomial $\binom{d-|M| + u}{u} = \frac{1}{u!} \prod_{\ell=1}^u (d-|M|+\ell)$, which vanishes when $d+1 \leq |M| \leq d+u$, so we may expand the condition in the sum to $|M| \leq d+u$. Now there are $n-s-d$ $\circ$'s and $s-u$ grey columns, so there are $n-(n-s-d)-(s-u) = d+u$ non-grey columns without $\circ$'s. Hence $|\Pi| \leq d+u$, so $M \subset \Pi$ automatically satisfies $|M| \leq d+u$ and we may remove the constraint on $|M|$ altogether. Since $|\Pi| > u$, the sum is thus zero by \Cref{cor:shifted_Vandermonde}.
    
    We claim that every term in \eqref{eq:hook_witness_extreme.2} has now been canceled precisely once. The ``orbits'' above obtained by applying the $\fS_s$-action to terms with $\Omega = \varnothing$ followed by relation (A) partition the terms, since starting at an arbitrary term, we may reverse the application of relation (A), which increases $d$ and therefore remains in $M_I$, to arrive at a term with $\Omega = \varnothing$. Terms in the orbit of $\Omega = \varnothing$ with $|\Pi| > u$, or equivalently terms with $|\Omega \sqcup \Pi| > u$, are thus entirely accounted for. For terms with $|\Omega \sqcup \Pi| = u$, we may first apply relation (A) to replace $\Omega$ with $\varnothing$ and $\Pi$ with $\Omega \sqcup \Pi$, then relation (B) as noted above applies to the shortest column with an $\circ$, so the same is true without needing to apply relation (A), resulting in a sign-reversing involution in the case $|\Omega \sqcup \Pi| = u$. Since $\Psi$ is a subset of the $s-u$ grey columns, we have $|\Omega \sqcup \Pi| = s - |\Psi| \geq s - (s-u) = u$, so all cases have been handled. This completes the proof.
  \end{proof}
\end{reptheorem}

\section{Further directions}\label{sec:further}

The lex-minimal $J$ appearing in either the Generic Pieri Rule, \Cref{thm:generic_Pieri}, or the extreme hook relations, \Cref{thm:hook_witness_extreme}, is $J=I$. Hence one may be tempted to use reverse lexicographic order on $2^{[n-1]}$ when attempting to answer \Cref{que:main}.

However, computations with $n=8$ show that this order together with the bijection $\Phi_n$ have correct composition factors at only $115$ out of $128$ cases. One may slightly tweak the reverse lexicographical order and get the predicted multiset of composition factors. For example, at $n=8, k=5$, replacing the reverse lex-interval
  \[ \{1, 2, 4, 5, 7\}, \{1, 2, 4, 5, 6\}, \{1, 2, 3, 6, 7\}, \{1, 2, 3, 5, 7\} \]
with
  \[ \{1, 2, 3, 6, 7\}, \{1, 2, 4, 5, 7\}, \{1, 2, 3, 5, 7\}, \{1, 2, 4, 5, 6\} \]
gives an affirmative answer to \Cref{que:main} in this case. In this way, orders verifying \Cref{que:main} valid for $n \leq 8$ have been found.

For $k \geq 2$, additional relations beyond those in our two families are required. A particular relation which is not explained by the results above is
\begin{align*}
0 &= 4\partial_{e_6(\underline{6})} \dif_{356} \Delta_8
-8\partial_{e_5(\underline{6})} \dif_{357} \Delta_8
+4\partial_{e_4(\underline{6})} \dif_{367} \Delta_8 \\
&-3\partial_{e_5(\underline{6})} \dif_{456} \Delta_8
+6\partial_{e_4(\underline{6})} \dif_{457} \Delta_8
-3\partial_{e_3(\underline{6})} \dif_{467} \Delta_8
\end{align*}

From our results and computations, the $\mathbb{Q}$-linear relations between $\partial_{e_r(\underline{m})} \dif_I \Delta_n$ exhibit rich combinatorial structure. Given the wealth of algebraic and geometric structure surrounding the various coinvariant algebras, we are led to the following.

\begin{problem}
  Completely describe the $\mathbb{Q}$-linear relations between $\partial_{e_r(\underline{m})} \dif_I \Delta_n$'s.
\end{problem}

\begin{problem}
  Give a conceptual explanation for the existence of these relations, perhaps in topological or geometric terms.
\end{problem}

\bibliography{refs}{}
\bibliographystyle{acm}

\end{document}